\magnification=\magstep1
\input amstex

\documentstyle{amsppt}

\voffset =-3pc

\NoBlackBoxes

\define\rr0{\operatorname{rr\;0}}

\define\sr1{\operatorname{sr\;1}}
\define\wtB{\widetilde{B}}
\define\wtC{\widetilde{C}}
\define\wtD{\widetilde{D}}
\define\wtJ{\widetilde{J}}
\define\wtI{\widetilde{I}}
\define\index{\operatorname{index}}

\define\tsr{\operatorname{tsr}}
\define\id{\operatorname{id}}
\define\her{\operatorname{her}}
\define\wtA{\widetilde{A}}

\define\bb1{\text{\bf 1}}

\document

\topmatter
\title
Non-Stable K--Theory and\\
Extremally Rich C$^*$-Algebras
\endtitle
\author
Lawrence G. Brown and Gert K. Pedersen
\endauthor
\abstract{We consider the properties
weak cancellation, $K_1-$surjectivity, good index theory, and
$K_1-$injectivity, for the class of extremally rich
$C^*-$algebras, and for the smaller class of isometrically rich
$C^*-$algebras.  We establish all four properties for  isometrically
rich $C^*-$algebras and for extremally
rich $C^*-$algebras that are either purely infinite or of
real rank zero, $K_1-$injectivity in the real rank zero case following 
from 
a prior result of H. Lin.  We also show that weak cancellation implies 
the other properties for extremally rich 
$C^*-$algebras  
and that the class of extremally rich $C^*-$algebras 
with weak cancellation is closed under extensions.  Moreover, we consider
analogous properties which replace the group $K_1(A)$ with the
extremal $K-$set $K_e(A)$ as well as two versions of $K_0-$surjectivity.}
\endabstract
\endtopmatter

\bigskip

\subhead{1. Introduction}\endsubhead

\medskip

In [{\bf 9}] we defined the concept of extremal richness. 
One of several equivalent criteria for the $C^*-$algebra $A$ to be
extremally rich is that the closed unit ball of $\wtA$, the unitization,
is the convex hull of $\Cal E(\wtA)$, the set of its extreme points.  Further
review of the concept is given in the next section.  A simple
$C^*-$algebra is extremally rich if and only if it is either 
of stable rank one or purely infinite, and a theme of our 
work has been that extremal richness is a generalization of the stable rank
one property which is suitable for infinite algebras.  Since much of
the success in the classification of simple $C^*-$algebras has been for
algebras that are either purely infinite or of stable rank one, it
seems worthwhile to study non-simple extremally rich $C^*-$algebras.

In [{\bf 17}] J. Cuntz defined purely
infinite simple C*-algebras and showed that they have many good
non-stable K-theoretic 
properties.  
And in [{\bf 35}] M. Rieffel, motivated by algebraic results of H. 
Bass [{\bf 4}], defined topological stable rank and showed
that $C^*-$algebras of (topological) stable rank one have similarly 
good properties.
We therefore investigated whether extremally rich $C^*-$algebras
also have the good properties.
Although we haven't proved that all
extremally rich algebras have good non-stable $K-$theoretic properties,
we have found large subclasses that do.  In particular, 
the summary Theorem 6.10
includes all four properties listed in the abstract for isometrically
rich $C^*-$algebras and for extremally
rich $C^*-$algebras which are either purely infinite (in the sense of 
E. Kirchberg and M. R\o rdam) or
of real rank zero.
All three cases of Theorem 6.10 cover purely infinite simple $C^*-$algebras.

A $C^*-$algebra $A$ has {\it weak cancellation} if whenever $p$ and
$q$ are projections in $A$ which generate the same (closed, two-sided)
ideal $I$ and have the same class in $K_0(I)$, then $p$ is
Murray-von Neumann equivalent to $q$ ($p\sim q$).  Of course,
$p\sim q$ implies that they generate the same ideal $I$ and that
$[p]=[q]$ in $K_0(I)$.  Nevertheless, it was observed by Rieffel,
c\.f\. [{\bf 5}, 6.5.1], that $C^*-$algebras of stable rank one
($\tsr(A)=1$)
satisfy a stronger property: If $[p]=[q]$ in $K_0(A)$, then $p\sim q$.
Of course this stronger property can't hold in infinite algebras.
Cuntz showed in [{\bf 17}] that if $A$ is purely infinite simple,
if $p$ and $q$ are both non-zero, and if
$[p]=[q]$ in $K_0(A)$, then $p\sim q$; and the concept of weak
cancellation was designed to specialize to this property in the
simple case.  More about weak cancellation and its history is given
in the next section.

A central question for our work on this paper was whether every 
extremally rich $C^*-$algebra has weak cancellation.  We
didn't find the answer, and some remarks about 
this question appear in \S 7.3.
Corollary 3.6 below states that every isometrically rich $C^*-$algebra has
weak cancellation.  Isometric richness, which is reviewed in the next
section, is a concept intermediate between stable rank one and extremal
richness.  It is equivalent to extremal richness for prime $C^*-$algebras,
in particular for simple ones.  Corollary 3.7 {\it a fortiori} proves
weak cancellation in the purely infinite case, and Theorem 3.10
is the real rank zero case.  Some closure properties of weak cancellation
within the class of extremally rich $C^*-$algebras, in particular
invariance under extensions, are summarized in Theorem 4.5.
Theorem 3.5 gives a structural property which is equivalent to
weak cancellation for extremally rich $C^*-$algebras, and the 
structural consequences of weak cancellation are further elaborated
in Proposition 3.15 and in 3.17--3.20.

We say that $A$ has $K_1-${\it surjectivity} if the map from
$\Cal U(\wtA)/\Cal U_0(\wtA)$ to $K_1(A)$ is surjective,
$K_1-${\it injectivity} if this map is injective, and $K_1-${\it bijectivity} if
it is bijective.  Here $\Cal U(\wtA)$ is the unitary group and 
$\Cal U_0(\wtA)$ the connected component of the identity (so that
$\Cal U(\wtA)/\Cal U_0(\wtA)$ is the set of homotopy classes).  Rieffel
showed in [{\bf 35}] that $\tsr(A)=1$ implies that $A$ has 
$K_1-$bijectivity and Cuntz [{\bf 17}] showed the same for $A$ purely
infinite simple.  A result of P. Ara, [{\bf 1}, Theorem 3.5]
states that every quotient $C^*$--algebra
of a Rickart $C^*$--algebra has $K_1$--surjectivity.
Theorem 4.4 below states that every extremally rich
$C^*-$algebra with weak cancellation has $K_1-$surjectivity.
This should be compared with Theorem 3.1 of [{\bf 3}], which states
that every $C^*-$algebra with real rank zero and stable weak
cancellation has $K_1-$surjectivity.

In [{\bf 11}, Definitions 3.6] we defined the extremal $K-$set, $K_e(A)$,
which is roughly analogous to $K_1(A)$ with extremal partial isometries used 
in place of unitaries.  The equivalence relation for $K_e$ is more
complicated than that for $K_1$; but in Corollary 5.3 we show that if
$A$ is extremally rich with weak cancellation, then 
$K_e(A)=\varinjlim_n(\Cal E(\Bbb M_n(\wtA))/\text{homotopy})$, in exact
analogy with the  $K_1-$case.
We say that $A$ has $K_e-${\it surjectivity}, $K_e-${\it injectivity}, 
or $K_e-${\it bijectivity} if the map from $\Cal E(\wtA)/\text{homotopy}$ to $K_e(A)$
is respectively surjective, injective, or bijective.  These properties
actually imply the corresponding $K_1-$properties.  In Theorem 4.7 we
show that every extremally rich 
$C^*-$algebra with weak cancellation has $K_e-$surjectivity.

We say that a (non-unital) $C^*-$algebra $K$ has {\it good index theory} if 
whenever $K$ is embedded as an ideal in a unital $C^*-$algebra $A$ and
$u$ is a unitary in $A/K$ such that $\partial_1([u]_{K_1})=0$ in
$K_0(K)$, there is a unitary in $A$ which lifts $u$.  Of course, the
boundary map, $\partial_1:K_1(A/K)\to K_0(K)$, from the $K-$theory long
exact sequence is often regarded as an abstract index map.  Using
that long exact sequence, we can reformulate the good index theory
property:  If $u\in \Cal U(A/K)$, $\alpha \in K_1(A)$, and $\alpha$ lifts
$[u]_{K_1}$, then $u$ lifts to $\Cal U(A)$.  This suggests a stronger
property--require that the lift of $u$ lie in the class $\alpha$.
However, there is no need to name this stronger property (still
demanded for all choices of $A$), since it is equivalent to requiring that
$K$ have both $K_1-$surjectivity and good index theory.

Good index theory has been considered by other authors, but, so far as
we know, no one previously proposed a name for it.  M. Pimsner, S. Popa,
and D. Voiculescu proved in [{\bf 34}, Lemma 7.5] that
$C(X)\otimes \Bbb K$ has good index theory when $X$ is compact, where
$\Bbb K$ denotes the algebra of compact operators on a separable,
infinite dimensional Hilbert space.  J. Mingo showed  in [{\bf 28}, Proposition 
1.11] that $C(X)$ can be replaced by an arbitrary unital $C^*-$algebra and
asked whether every stable $C^*-$algebra has good index theory.  (In both
of these results the property is considered only for $A=M(K)$, the 
multiplier algebra, but this is sufficient to imply that it holds for
all $A$.)  Shortly after, G. Nagy proved ([{\bf 29}, Theorem 2], 
{\it a fortiori}), that $\text{csr}(K)\le 2$ implies that $K$ has
good index theory.  Here csr denotes connected stable rank (Rieffel [{\bf 35}]),
and $\text{csr}(K)\le 2$ also implies $K_1-$surjectivity for $K$.
  It had already been proved by A. Sheu, [{\bf 38}, Theorem 3.10], and
V. Nistor, [{\bf 30}, Corollary 2.5], that $\text{csr}(K)\le 2$ for all
stable $K$.  So Nagy completed the affirmative answer to Mingo's question 
(but was apparently unaware of Mingo's paper).  Since Rieffel proved in
[{\bf 35}] that $\text{csr}(K)\le \tsr(K)+1$, Nagy also established good
index theory for $C^*-$algebras of stable rank one.  The fact that purely
infinite simple $C^*-$algebras have good index theory should also be
considered previously known, but we haven't found a precise reference
for it (c\.f\. [{\bf 43}, \S2], especially the proof of Lemma 2.2).

Theorem 5.1 below states that every extremally rich $C^*-$algebra with
weak cancellation has good index theory.  (An analogous result for 
$C^*-$algebras of real rank zero was independently discovered by the first 
named author and 
F. Perera, cf\. [{\bf 33}, Theorem 3.1].)  We also consider several
analogues of good index theory which use the extremal $K-$set in place of
$K_1$.  These are spelled out, but not named, in 5.4, and several 
results and remarks about them are in \S5.

The main result of \S6 is that every extremally rich 
$C^*-$algebra with weak cancellation has $K_1-$injectivity, which is in
Theorem 6.7.  S. Zhang proved in [{\bf 41}] (an alternative to the 
actually earlier proof in [{\bf 42}, Theorem 1.2]) that purely infinite 
simple implies real rank zero, and H. Lin proved in [{\bf 27}, Lemma 2.2] 
that real rank zero implies $K_1-$injectivity.  Thus Lin's result, as 
well as ours, is a generalization of Cuntz's $K_1-$injectivity result in
[{\bf 17}].  Our proof requires the use of (strong and weak) 
$K_0-$surjectivity, which are defined in 6.2 below.  We also show in \S6 
that every extremally rich 
$C^*-$algebra with weak cancellation has weak $K_0-$surjectivity and 
$K_e-$injectivity.

This is the last of our joint papers.  It is the sixth paper of a project
that began in 1992, the others being [{\bf 9}], [{\bf 10}], [{\bf 11}], 
[{\bf 12}], and [{\bf 13}].
Initially there were to be only two papers.  The first paper was split
into [{\bf 9}], which consists of its first six sections, and [{\bf 10}], which
consists of the remaining three sections and an introduction.  
This paper,
[{\bf 12}], and [{\bf 13}] constitute an expanded version of the original second paper, 
and [{\bf 11}] is based on an idea that came later.  Several of the main results
of the present paper, Corollary 3.6 and Theorems 4.4 and 5.1, were obtained
in 1992, but the proofs have been improved since then.

We are grateful to P. Ara, K. Goodearl, and F. Perera for helpful comments.

\subhead{2. Notations and Preliminaries}\endsubhead

\medskip

\definition{2.1. Extremal Richness, etc} A unital $C^*-$algebra $A$ has 
{\it stable rank one} if $A^{-1}$, the set of invertible elements, is 
dense in $A$, is {\it isometrically rich} if $A_l^{-1}\bigcup A_r^{-1}$, 
the set of one-sided invertible elements, is dense, and is 
{\it extremally rich} if $A_q^{-1}$, the set of quasi-invertible 
elements is dense.  If $A$ is non-unital, we say that $A$ has one of 
these properties if $\wtA$ has the property.
All three properties pass to (closed, two-sided) ideals and hereditary $C^*-$subalgebras, 
quotient algebras, and matrix algebras and stabilizations.

Before reviewing the definition of quasi-invertibility, 
we recall R. Kadison's criterion ([{\bf 22}]) for extreme points of the
unit ball of a $C^*-$algebra $A$.  Extreme points exist if and only if 
$A$ is unital, and $u$ is extremal if and only if 
 $$
(\bb1-uu^*)A(\bb1-u^*u)=0.
 $$
Equivalently, $u$ is a partial isometry and $I\cap J=0$, where $I=\id(\bb1-uu^*)$
and $J=\id(\bb1-u^*u)$.  Here $\id(\cdot)$ denotes the ideal 
generated by $\cdot$, $\bb1-uu^*$ and $\bb1-u^*u$ are called the left 
and right {\it defect projections} of $u$, and $I$ and $J$ are the 
left and right {\it defect ideals} of $u$.

In [{\bf 9}, Theorem 1.1] we showed that seven conditions on an 
element $t$ of a unital $C^*-$algebra $A$ are equivalent, and these 
conditions are the definition of {\it quasi-invertible}.  (Non-unital 
algebras have no quasi-invertibles.)  One of these conditions amounts
to saying that $t$ has closed range and that if $t=u|t|$ is its canonical
polar decomposition (the closed range condition implies $u\in A$), then 
$u\in \Cal E(A)$.
Another is that 
there  are ideals $I$ and $J$ with $I\cap J=0$ such that $t+J$ is left 
invertible in $A/J$ and $t+I$ is right invertible in $A/I$.
Clearly the minimal choices for $I$ and $J$ are the defect ideals of $u$, 
and we also call these the {\it defect ideals} of $t$. 
Of course, if $A$ is prime, one of $I$ and $J$ must be $0$; thus every 
quasi-invertible element is one-sided invertible and every extremal 
partial isometry is an isometry or co-isometry.

In general, there is an analogy in the two-step progressions from stable rank
one (through isometrically rich) to extremally rich, from invertibility 
to quasi-invertibility, and from unitary to extremal partial isometry.
It is not always necessary to pursue all three concepts in parallel
because $A$ (unital) has stable rank one if and only if it is extremally rich and 
every extremal is unitary, and $A$ is isometrically rich if and only if 
it is extremally rich and every extremal is an isometry or co-isometry 
(c\.f\. [{\bf 12}, Proposition 4.2]).  
However, the most general results are usually not proved first.
In [{\bf 36}, Corollary 3.7] R\o rdam proved Robertson's conjecture: 
For $A$ unital, $\tsr(A)=1$ if and only if the closed unit ball is
the convex hull of $\Cal U(A)$.
The second named author analyzed the situation with $\Cal U(A)$ replaced 
by $\Cal E(A)$ when $A$ is prime in [{\bf 32}, \S8], and our generalization
for arbitrary $A$ is in [{\bf 10}, \S3].  

Our contention that quasi-invertibility is a very natural concept can 
be bolstered by some other results from our earlier papers,
for example the easy Proposition 1.1 in [{\bf 12}] or the more
technical treatment of elements with persistently closed range in 
[{\bf 11}, \S 7].  Quasi-invertibility also plays a key role in an 
extension of the classical index theory of semi-Fredholm operators. 
This was given in [{\bf 11}, \S 6, 7]  and 
is briefly described in Remark 5.9(i) below.  It is a classical result that 
an element of $B(H)$ is semi-Fredholm if and only if it has 
persistently closed range.

For technical reasons it is sometimes necessary to consider extremals, 
quasi-invertibility, and extremal richness for objects other than 
$C^*-$algebras, namely bimodules of the form $pAq$, where $p$ and $q$ 
are projections in $A$.  S. Sakai's criterion for $u$ to be an extreme point of
the unit ball of $pAq$ (c\.f\. [{\bf 31}, 1.4.8]) is: 
 $$
(p-uu^*)A(q-u^*u)=0\text{ or equivalently, }(\bb1-uu^*)pAq(\bb1-u^*u)=0.
 $$
Quasi-invertibility and extremal richness for bimodules are treated 
analogously to the treatments for $C^*-$algebras with no difficulty.  It 
is not necessary to be explicitly aware of the fact that $pAq$ is a 
bimodule or even to know what ``bimodule'' means.  Nevertheless, 
the abstract setting was discussed in [{\bf 9}, \S4].  One warning: 
There are no concepts of unitality or unitization for bimodules.  It is
required that $\Cal E(pAq) \ne \emptyset$ in order for $pAq$ to have a
chance to be extremally rich.  (When $pAq=0$, it is automatically extremally
rich.)  The extremal richness of $A$ does not imply that of $pAq$, but it 
is important to our main results that sometimes $pAq$ is extremally rich 
(c\.f\. Proposition 3.2 below).  Since $\tsr(A)\le n$ if and only if 
``left invertibles'' are dense in the bimodule ${\wtA}^n$ 
($=\bb1_n\Bbb M_n(\wtA)\bb1_1$), we once looked at extremal richness for
${\wtA}^n$; but it turned out that $A$ extremally rich does not imply 
${\wtA}^n$ extremally rich.
\enddefinition
\medskip
\example{Aside}  An example that does invoke the abstract 
setting may be interesting.  The right module called $\Cal H_A$ by G. Kasparov 
in [{\bf 23}] is also an $A\otimes\Bbb K$--$A$--imprimitivity bimodule.  
Regardless of whether or not $A$ is extremally rich, 
$\Cal H_A$ is an extremally rich bimodule if and only if $A$ is unital.
\endexample

\medskip

\definition{2.2. Notations and Definitions} If $p$ and $q$ are projections 
in a $C^*-$algebra, we write $p\precsim q$ to mean $p\sim q'\le q$ for
some projection $q'$.  
The relations $\sim$ and $\precsim$  can be extended to projections in 
$\bigcup_n\Bbb M_n(A)$ either by replacing $p$ and $q$ by $p\oplus 0_k$ and 
$q\oplus 0_l$ for suitable $k$ and $l$ or by allowing the partial isometries 
to be non-square matrices. 
The projection $p$ is {\it infinite} if it is 
equivalent to a proper subprojection of itself, otherwise {\it finite}, 
and a unital $C^*-$algebra $A$ is finite or infinite according as $\bb1_A$ is,  
and {\it stably finite} if all the matrix algebras $\Bbb M_n(A)$ are finite.
If $p$ and $q$ are projections in $\Bbb M_m(A)$ and $\Bbb M_n(A)$, 
respectively, then $p\oplus q$ denotes the projection 
$\left(\smallmatrix p&0\\
0&q\endsmallmatrix\right)$ 
in $\Bbb M_{m+n}(A)$, and $p_1\oplus \dots \oplus p_k$ is defined similarly. 
In this context $kp$ is used for the $k-$fold sum, $p\oplus\dots\oplus p$. 
A projection $p$ is called {\it properly infinite} if $2p\precsim p$.
\enddefinition
\medskip

\definition{2.3. More on Weak Cancellation} 
It is possible to reformulate weak cancellation in a way that does not
mention $K-$theory.  Note that if $p$, $q$, and $r$ are projections in 
an ideal $I$  and if $p$ is full in $I$ (i\.e\., $I=\id(p)$), then
$[q]_{K_0(I)}=[r]_{K_0(I)}$ if and only if $q\oplus np\sim r\oplus np$ for
sufficiently large $n$.  Thus the hypotheses, $\id(p)=\id(q)=I$ and 
$[p]_{K_0(I)}=[q]_{K_0(I)}$, can be replaced by, $p\oplus nq\sim (n+1)q$ 
and $q\oplus np\sim (n+1)p$ for sufficiently large $n$.  It was pointed out 
to us by K. Goodearl that the concept can be simplified further if we 
demand weak cancellation for the stabilization of $A$:  $A\otimes \Bbb K$ has 
weak cancellation if and only if $2p\sim p\oplus q\sim 2q$ implies
$p\sim q$ for all projections $p$ and $q$ in $A\otimes \Bbb K$.  Moreover, 
this is 
equivalent to ``separativity,'' a term which was introduced into 
semigroup theory by A. Clifford and G. Preston [{\bf 15}].  
The set of Murray-von Neumann equivalence classes of 
of projections in $A$ is not in general a semigroup, and it is only the
stable version of weak cancellation that is literally equivalent to 
separativity (c\.f\. 3.1 below).  However, [{\bf 2}, Theorem 2.8] can be 
used to show that weak cancellation is a stable property in the real rank
zero case.
More detailed
discussion of the history of separativity can be found in the introduction of
[{\bf 2}].
\enddefinition
\medskip

\definition{2.4. Defect Ideals}  Defect ideals were treated in 
[{\bf 12}, \S3].  The {\it defect ideal} of $A$, denoted $\Cal D(A)$, 
is the ideal generated by all defect projections of elements of 
$\Cal E(\wtA)$.  (All of these defect projections are in $A$.)  
If $A$ is extremally rich, then $\Cal D(A)$ is the 
smallest ideal $I$ such that $\tsr(A/I)=1$.  It follows that defect ideals
are compatible with Rieffel-Morita equivalence among extremally rich 
$C^*-$algebras.  In particular, if $B$ is a full hereditary 
$C^*-$subalgebra, then $\Cal D(B)=B\cap\Cal D(A)$.  The symbol 
$\Cal D^n(A)$ denotes the $n-$fold iteration of $\Cal D$.
\enddefinition
\medskip

\definition{2.5. Notations}  The primitive ideal space of $A$ will be 
denoted by $A^{\vee}$.  If $I$ is an ideal of $A$, then $I^{\vee}$ 
is identified with an open subset of $A^{\vee}$ (namely, the complement of 
$\text{hull}(I)$), and if $B$ is a hereditary
$C^*-$subalgebra, then $B^{\vee}$ is identified with $\id(B)^{\vee}$.
For $B$ hereditary, $B^{\perp}$ denotes the two-sided annihilator of
$B$, which is again hereditary, and $I^{\perp}$ is an ideal if $I$ is
an ideal.

Also, ${}^=$ denotes norm closure and $\Cal T_e$ denotes the extended 
Toeplitz algebra, which was discussed in [{\bf 9}, page 143].
\enddefinition

\subhead{3. Weak Cancellation}\endsubhead

\medskip

\definition{3.1. Definitions} Recall that a $C^*$--algebra $A$ has 
weak cancellation if any
pair of projections $p,q$ in $A$ that generate the same closed ideal
$I$ of $A$ and have the same image in $K_0(I)$ must be Murray--von
Neumann equivalent in $A$ (hence in $I$). If $\Bbb M_n(A)$ has weak
cancellation for every $n$, equivalently, if $A\otimes\Bbb K$ has
weak cancellation, we say that $A$ has {\it stable weak
cancellation}.  We shall show below that weak cancellation
implies stable weak cancellation if $A$ is extremally rich, but for
now we need the distinction.
\enddefinition
\bigskip

\proclaim{3.2. Proposition} If $p$ and $q$ are projections in an
extremally rich $C^*$--algebra $A$ such that $[p]=[q]$ in $K_0(A)$
then $pAq$ is extremally rich.
\endproclaim

\demo{Proof} Since $K_0(A)\subset K_0(\wtA)$ and $pAq=p\wtA q$ we may
assume that $A$ is unital. Then also $[\bb1-p]=[\bb1-q]$ in $K_0(A)$,
so
 $$
(\bb1-p)\oplus n\bb1\sim(\bb1-q)\oplus n\bb1 \quad (\text{in} \quad 
\Bbb M_{n+1}(A))
 $$
for $n$ sufficiently large. Since $\Bbb M_{n+1}(A)$ is extremally
rich by [{\bf 9}, Theorem 4.5] we can use [{\bf 9}, Proposition 4.2]
to conclude that
 $$
pAq=(p\oplus 0)\Bbb M_{n+1}(A)(q\oplus 0)
 $$
is extremally rich. \hfill$\square$
\enddemo

\bigskip

\proclaim{3.3. Lemma} Let $p$ and $q$ be projections in a
$C^*$--algebra $A$ and for each element $x$ in $A$ let $id(x)$ denote
the closed ideal generated by $x$. If now $v\in\Cal E(pAq)$ then
 $$
id(v)=id(p)\cap id(q)\;.
 $$
\endproclaim

\demo{Proof} Since $vv^*\le p$ we have $vv^*\in id(p)$, whence $v\in
id(p)$. Similarly $v^*v\le q$, so $v\in id(q)$.

Conversely, let $\pi:A\to A/id(v)$ denote the quotient map. Then the
extremality equation gives
 $$
\pi(p)\pi(A)\pi(q)=0\;,
 $$
whence $\pi(id(p)\cap id(q))=0$; so $id(p)\cap id(q)\subset id(v)$, as
desired. \hfill$\square$
\enddemo

\bigskip

Although we will prove later that every extremally rich $C^*-$algebra
with weak cancellation also has
$K_1$--surjectivity, 
it will facilitate some of the following
arguments to impose it as a condition on the algebras. 
\bigskip

\proclaim{3.4. Lemma} Let $p$ and $q$ be full projections in an
extremally rich $C^*$--algebra $A$ such that $[p]=[q]$ in $K_0(A)$,
and assume that we have found an extreme partial isometry $u$ in
$\Cal E(pAq)$ such that
 $$
p_1=p-uu^*\quad\text{and}\quad q_1=q-u^*u
 $$
are projections in an ideal $I$ of $A$. Assume further that
$eAe/eIe$ has $K_1$--surjectivity for every full projection $e$ in
$A$. Then
 $$
p\sim p_2\oplus e_2\quad\text{and}\quad q=q_2+e_2
 $$
for some full projection $e_2$ in $A$ and projections $p_2$ and
$q_2$ in $I$, with $[p_2]=[q_2]$ in $K_0(I)$.
 \endproclaim

\demo{Proof} Note first that the element $[p_1]-[q_1]$ of $K_0(I)$
belongs to the kernel of the natural map from $K_0(I)$ into $K_0(A)$.
By the six--term exact sequence in $K$--theory there is therefore an
$\alpha$ in $K_1(A/I)$ such that $\partial_1\alpha=[p_1]-[q_1]$. Since
$e_1=u^*u$ is a full projection in $A$ by Lemma 3.3, we may identify
$K_1(A/I)$ with $K_1(e_1Ae_1/e_1Ie_1)$, and by $K_1$--surjectivity we
can therefore find $v$ in $e_1Ae_1$ such that $v+e_1Ie_1$ is unitary
with $[v+e_1Ie_1]=\alpha$. Since $e_1Ae_1$ is extremally rich,
extreme points lift from quotients (cf\. [{\bf 9}, Theorem 6.1]), so
we may assume that $v\in\Cal E(e_1Ae_1)$. Computing in $K_0(I)$ we
find that
 $$
[p_1]-[q_1]=\partial_1\alpha=\index v=[e_1-v^*v]-[e_1-vv^*]\;.
 $$
Let
 $$
p_2=p_1+u(e_1-vv^*)u^*\;,\;\;q_2=q_1+(e_1-v^*v)\;,\;\;e_2=v^*v\;.
 $$
Since $v\in\Cal E(e_1Ae_1)$ we see from Lemma 3.3 that $e_2$ is full
in $e_1Ae_1$, and since $e_1$ is full in $A$, we have also that
$e_2$ is a full projection in $A$. By construction $e_1-vv^*$ and
$e_1-v^*v$ belong to $I$, so $p_2$ and $q_2$ belong to $I$; and
evidently $[p_2]=[q_2]$ in $K_0(I)$. Finally,
 $$
\aligned
p & = p_1+uu^*=p_2+uvv^*u^*\sim p_2\oplus e_2\;,\\
q & = q_1+u^*u=q_2+e_2\;.
\endaligned
 $$
\hfill$\square$
\enddemo

\medskip

\proclaim{3.5. Theorem} For an extremally rich $C^*$--algebra $A$
the following conditions are equivalent:
\roster
 \item"(i)" $A$ has weak cancellation;
 \item"(ii)" If $B=pAp$ for some projection $p$ in $A$ and $u\in\Cal
E(\Bbb M_2(B))$, there is a projection $q$ in $B$ such that
 $$
q\oplus 0\sim (p\oplus p) - uu^*\quad (\text{in} \quad \Bbb M_2(B))\;;
 $$
 \item"(iii)" If $B=pAp$ for some projection $p$ in $A$, and
$\{u_1,\dots,u_n\}$ is a finite subset of $\Cal E(B)$ there is a
projection $q$ in $B$ such that
 $$
q\oplus 0\sim\bigoplus_{i=1}^n(p-u_iu_i^*)\quad (\text{in}\quad
\Bbb M_n(B)) \;.
 $$
\endroster
\endproclaim

\demo{Proof} (i) $\Rightarrow$ (ii) Since homotopic projections are
equivalent we may assume by \cite{{\bf 11}, Corollary 2.13} that
 $$
(p\oplus p) - uu^* = (p-v_1v_1^*-w_{12}w_{12}^*)\oplus
(p-v_2v_2^*-w_{21}w_{21}^*)\;,
 $$
where $v_i\in\Cal E(B)$ and $w_{ij}\in\Cal E(p_iBq_j)$, and where we
set $p_i=p-v_iv_i^*$ and $q_i=p-v_i^*v_i$ for $i,j=1,2$. Evidently
the two projections $p-w_{21}w_{21}^*$ and $p-w_{21}^*w_{21}$ have
the same image in $K_0(I)$ for any ideal $I$ containing $p$. To show
that they also generate the same ideal in $A$ it suffices to
show that they both generate $B$ as an ideal (inside $B$). However,
 $$
p-w_{21}w_{21}^*\ge v_2v_2^*\quad\text{and}\quad
p-w_{21}^*w_{21}\ge v_1^*v_1\;,
 $$
and both $v_1$ and $v_2$ generate $B$ as an ideal by Lemma 3.3.

By assumption there is therefore a partial isometry $v$ in $B$ such
that
 $$
v^*v=p-w_{21}w_{21}^*\quad\text{and}\quad
vv^*=p-w_{21}^*w_{21}\;.
 $$
Thus
 $$
e=v(p_2-w_{21}w_{21}^*)v^*
 $$
is a projection equivalent to $p_2-w_{21}w_{21}^*$, and since $e\le
vv^*=p-w_{21}^*w_{21}$ and $e$ is centrally orthogonal to
$q_1-w_{21}^*w_{21}$ (because $w_{21}\in\Cal E(p_2Bq_1)$) it follows
that actually
 $$
e\le p-q_1=v_1^*v_1\;.
 $$
But then $v_1ev_1^*$ is a projection in $B$ equivalent to
$p_2-w_{21}w_{21}^*$ and orthogonal to $p-v_1v_1^*=p_1$. We may
therefore take
 $$
q=p_1-w_{12}w_{12}^*+v_1ev_1^*\;,
 $$
which is a projection in $B$ equivalent to $(p\oplus p) - uu^*$.

(ii) $\Rightarrow$ (iii) We use induction on $n$, the case $n=1$
being trivial. By assumption, used on $\{u_1^2,\dots,u_{n-1}^2\}$,
we can therefore find a projection $q_1$ in $B$ such that
 $$
q_1\sim\bigoplus_{i=1}^{n-1}(p-u_i^2u_i^{*2})\sim
2\bigoplus_{i=1}^{n-1}(p-u_iu_i^*)\;.
 $$
This means that we can write $q_1=q_2+(q_1-q_2)$, where
 $$
q_2\sim q_1-q_2\sim\bigoplus_{i=1}^{n-1}(p-u_iu_i^*)\;.
 $$
In particular, $q_2\precsim p-q_2=p_1$.

Set $B_1=p_1Ap_1$. Since
 $$
p=p_1+q_2\precsim 2p_1
 $$
there is a partial isometry $w$ in $\Bbb M_2(B)$ such that $w^*w=p\oplus 0$
and $ww^*=e\le p_1\oplus p_1$. It follows that if we define
 $$
v=p_1\oplus p_1-e+w(u_n\oplus 0)w^*
 $$
then $v\in\Cal E(\Bbb M_2(B_1))$. Moreover,
 $$
p_1\oplus p_1-vv^*=w((p-u_nu_n^*)\oplus 0)w^*\sim p-u_nu_n^*\;.
 $$
 
Applying (ii) to $B_1$ and $v$ we find a projection $p_2$ in $B_1$
such that $p_2\sim p-u_nu_n^*$. Thus we may take $q=q_2+p_2$ to
complete the induction step.

(iii) $\Rightarrow$ (i) Let $p$ and $q$ be projections in $A$ that
generate the same closed ideal $I$ and for which $[p]=[q]$ in
$K_0(I)$. Replacing $A$ by $I$ does not effect the conditions in
(iii) so we may assume that $I=A$; i\.e\., $p$ and $q$ are full
projections in $A$. Since $[p]=[q]$ in $K_0(A)$ it follows from
Proposition 3.2 that $pAq$ is extremally rich (and non--zero). Take
therefore $u$ in $\Cal E(pAq)$ and define
 $$
p_1=p-uu^*\;,\quad q_1=p-u^*u\;.
 $$
 Then $[p_1]=[q_1]$ in $K_0(A)$. 

If $\pi:A\to
A/\Cal D(A)$ denotes the quotient map, we have $[\pi(p)]=[\pi(q)]$ in
$K_0(\pi(A))$. Since $\tsr(\pi(A))=1$, this
implies that $\pi(p)\sim\pi (q)$ by [{\bf 5}, 6.5.1]. Thus $\pi(pAq)$
is isometrically isomorphic to $\pi(pAp)$, which has stable rank one.
It follows that $\pi(u)$ is ``unitary'' so that
$\pi(p_1)=\pi(q_1)=0$. Thus both $p_1$ and $q_1$ belong to $\Cal
D(A)$. 

Since $\tsr(eAe/e\Cal D(A)e)=1$ for every projection $e$ in $A$ we
may apply Lemma 3.4 with $I=\Cal D(A)$ to obtain a full projection
$e_2$ in $A$ and projections $p_2$ and $q_2$ in $\Cal D(A)$ such that
$[p_2]=[q_2]$ in $K_0(\Cal D(A))$ and
 $$
p\sim p_2\oplus e_2\;,\quad q=q_2+e_2\;.
 $$
 
Let $B=e_2Ae_2$. Then $B$ is a full corner of $A$, so $e_2\Cal
D(A)e_2=\Cal D(B)$ is a full hereditary $C^*$--subalgebra of $\Cal
D(A)$. Consequently the set of projections
 $$
\Cal D=\{e_2-ww^*\mid w\in\Cal E(B)\}
 $$
generates $\Cal D(A)$ as an ideal. For some finite subset $\{w_i\}$
of $\Cal D$ we therefore have
 $$
p_2\precsim\bigoplus(e_2-w_iw_i^*)\;,\;\;
q_2\precsim\bigoplus(e_2-w_iw_i^*)\;;
 $$
and since $[p_2]=[q_2]$ in $K_0(\Cal D(A))$ we can assume, possibly
after enlarging the subset, that
 $$
p_2\oplus(\bigoplus(e_2-w_iw_i^*))\sim q_2\oplus(\bigoplus(e_2-w_iw_i^*))\;.
 $$
Applying condition (iii) we find a projection $q_0$ in $B$ such that
$q_0\sim\bigoplus(e_2-w_iw_i^*)$. Thus $p_2\oplus q_0\sim q_2+q_0$ with
$q_0\le e_2$, whence
 $$
p\sim p_2\oplus e_2\sim q_2+e_2=q\;.
 $$
\hfill$\square$
\enddemo

\medskip

\proclaim{3.6. Corollary} Every isometrically rich $C^*$--algebra
has stable weak cancellation.
\endproclaim

\demo{Proof} If $A$ is isometrically rich then so is $\Bbb M_n(A)$
for every $n$ by \cite{{\bf 12}, Proposition 4.5} in conjunction with [{\bf 9}, Theorem
4.5]. It suffices therefore to show that $A$ has weak cancellation.
But this is evident since condition (iii) in our previous theorem is
covered by \cite{{\bf 12}, Lemma 4.6}. \hfill$\square$
\enddemo

\medskip
The term ``purely infinite'' will be used in the Kirchberg-R\o rdam
sense, c\.f\. \cite{\bf 24}, \cite{\bf 25}.  Every projection in
a purely infinite $C^*-$algebra is properly infinite.
In \cite{{\bf 12}, 3.10} we made the definition that $A$ is
{\it purely properly infinite} if every non-zero hereditary
$C^*-$subalgebra is generated as an ideal by its properly infinite
projections.  This is equivalent to purely infinite 
when $A$ is simple but stronger  in general. 
However, the two concepts are equivalent if $A$ has enough projections.  In \cite{{\bf 12}, Theorem 3.9}
we showed that an extremally rich $C^*-$algebra $A$ is purely properly
infinite if and only if $\Cal D(I)=I$ for every ideal $I$ of $A$.  We
also showed in \cite{{\bf 12}, Lemma 3.8} that for $A$ extremally rich
every non-zero projection in $A$ is properly infinite if and only if
$\Cal D(I)=I$ for every ideal $I$ which is generated (as an ideal) by
a projection.  
The next result has a still weaker hypothesis.
\bigskip

\proclaim{3.7. Corollary} If $A$ is an extremally rich $C^*$--algebra 
such that $\Cal D(I)=I$ whenever $I$ is the left defect ideal of an
element of $\Cal E(\wtA)$, then $A$ has weak cancellation. In particular,
this applies if every defect projection is properly infinite or
if $A$ is purely infinite. 
\endproclaim

\demo{Proof} 
To show that $A$ satisfies condition (iii) let
$u_1,\dots,u_n$ be in $\Cal E(B)$, where $B=pAp$ for a projection $p$
in $A$, and let $q_i=p-u_iu_i^*$.  Since $q_i$ is also a defect projection
for $A$, the hypothesis implies that $\Cal D(\text{id}(q_i))=\text{id}(q_i)$;
and since $q_iBq_i$ is Rieffel-Morita equivalent to $\text{id}(q_i)$,
this implies $\Cal D(q_iBq_i) =q_iBq_i$.  Then [{\bf 12}, Lemma 3.5]
implies that $mq_i$ is properly infinite for some $m$.  But $mq_i$
is equivalent to a projection $r_i$ in $B$, namely the left defect
projection of $u_i^m$.  Now the proper infiniteness of $r_i$ implies
that there is an isometry $v_i$ in $r_iBr_i$ such that
$r_i\precsim r_i-v_iv_i^*$.  So if $w_i=p-r_i+v_i$, then $w_i$ is an
isometry in $B$ and $q_i\precsim p-w_iw_i^*$.  Finally, if
$w=\prod_{i=1}^nw_i$, then $w$ is an isometry in $B$ and
$\oplus_{i=1}^nq_i\precsim (p-ww^*)$.
\hfill$\square$
\enddemo

\medskip

\proclaim{3.8. Lemma} Let $p,q$ and $q_0$ be projections in a unital
$C^*$--algebra $A$ such that
 $$
q\sim q_0\le p\quad\text{and}\quad
\bb1-q\sim\bb1-p\sim\bb1\;.
 $$
Then also $\bb1-q_0\sim\bb1$.
\endproclaim

\demo{Proof} We compute (in $A$ and in $\Bbb M_2(A))$
 $$
\aligned
& \bb1-q_0=\bb1-p+p-q_0\sim(\bb1-p)\oplus(p-q_0)\\
& \sim\bb1\oplus(p-q_0)=(\bb1-q+q)\oplus(p-q_0)\\
& \sim(\bb1-p+q_0)\oplus p-q_0\sim\bb1-p+q_0+p-q_0=\bb1\;.
\endaligned
 $$
\hfill$\square$ 
\enddemo

\medskip

\proclaim{3.9. Lemma} Let $q$ and $q_0$ be projections in a unital
$C^*$--algebra $A$ of real rank zero such that $q\sim q_0$ and
$\bb1-q\sim\bb1-q_0\sim\bb1$. If $I$ denotes the closed ideal of $A$
generated by $q$ (and $q_0$) then $\bb1-q\sim\bb1-q_0$ in
$\wtI=I+\Bbb C\bb1$. 
\endproclaim

\demo{Proof} By assumption there are $u$ and $v$
in $A$ such that
 $$
u^*u=v^*v=\bb1\;,\quad uu^*=\bb1-q\;,\quad vv^*=\bb1-q_0\;.
 $$
Let $\pi:A\to A/I$ be the quotient map and note that $\pi(u)$ and
$\pi(v)$ are unitaries. If $w$ is the partial isometry in $A$ for
which $w^*w=q_0$ and $ww^*=q$ (so that $w\in I$) then $uv^*+w$ is
unitary in $A$.

Since $\pi (A)$ has $K_1-$injectivity by Lin, [{\bf 27}, Lemma 2.2], 
and since $[\pi(v^*uvu^*)]=0$ in
$K_1(\pi(A))$, we must have
 $$
\pi(v^*uvu^*)=\pi(u_0)
 $$
for some unitary $\pi(u_0)$ in the identity component of $\pi(\Cal
U(A))$. However, $\Cal U_0(\pi(A))=\pi(\Cal U_0(A))$ so we may
assume that $u_0\in\Cal U_0(A)$. Consider now the unitary
$w_1=u_0(uv^*+w)$ and note that
 $$
\pi(w_1)=\pi(v^*uvu^*uv^*)=\pi(v^*u)\;.
 $$
Therefore $v_1=vw_1u^*$ is a partial isometry in $A$, and by
computation
 $$
v_1^*v_1=uu^*=\bb1-q\quad\text{and}\quad v_1v_1^*=vv^*=\bb1-q_0\;.
 $$
Finally, 
 $$
\pi(v_1)=\pi(v(v^*u)u^*)=\pi(\bb1)\;,
 $$
so that $v_1\in\wtI$, as desired. \hfill$\square$
\enddemo

\bigskip

\proclaim{3.10. Theorem} Every extremally rich $C^*$--algebra $A$ of
real rank zero has stable weak cancellation.
\endproclaim

\demo{Proof} The given data are stable so it suffices to show that
$A$ has weak cancellation. We do this by verifying condition (ii) in
Theorem 3.5, and since the given data are also hereditary it
suffices to verify the condition for $A$ alone, assuming that $A$ is
unital. Finally, using \cite{{\bf 11}, Corollary 2.13} we may assume that the defect
projection in $\Bbb M_2(A)$ has the form
 $$
(\bb1-v_1v_1^*-w_{12}w_{12}^*)\oplus
(\bb1-v_2v_2^*-w_{21}w_{21}^*)\;,
 $$
where $v_i\in\Cal E(A)$ and $w_{ij}\in\Cal E(p_iAq_j)$, for
$p_i=\bb1-v_iv_i^*$ and $q_i=\bb1-v_i^*v_i$, $i,j=1,2$.

Let $J$ be the closed ideal of $A$ generated by the two interesting
projections $p_1-w_{12}w_{12}^*$ and $p_2-w_{21}w_{21}^*$. We have
 $$
(p_2-w_{21}w_{21}^*)A(q_1-w_{21}^*w_{21})=0
 $$
since $w_{21}\in\Cal E(p_2Aq_1)$. Moreover,
 $$
(p_1-w_{12}w_{12}^*)A(q_1-w_{21}^*w_{21})=0
 $$
since already $p_1Aq_1=0$. It follows that $q_1-w_{21}^*w_{21}\in
J^\bot$. Since $J$ is isomorphic to its image in $A/J^\bot$ we may
replace $A$ by $A/J^\bot$ without changing notation. In other words
we may assume that $J^\bot=0$. In that case $q_1-w_{21}^*w_{21}=0$,
so if we put $q_0=w_{21}w_{21}^*$ we have $q_1\sim q_0\le p_2$.

Let $I$ be the closed ideal of $A$ generated by $q_1$ (and $q_0$).
Since $p_1Aq_1=0$ we see that $p_1\in I^\bot$. But since $q_1\precsim
p_2$ and $q_2Ap_2=0$ also $q_2\in I^\bot$. With $\pi:A\to A/I^\bot$
the quotient map this means that the three projections
 $$
\pi(p_2)\;,\quad\pi(q_1)\;,\quad\pi(q_0)
 $$
satisfy the assumptions of Lemma 3.8. Consequently
$\pi(\bb1-q_0)\sim\pi(\bb1)\sim\pi(\bb1-q_1)$. But now Lemma 3.9
applies to show that $\pi(\bb1-q_1)=u^*u$ and $\pi(\bb1-q_0)=uu^*$
for some partial isometry $u$ in $\pi(\wtI)$. Since $\pi(\wtI)$ is
isomorphic to $\wtI$ this means that $\bb1-q_1\sim\bb1-q_0$ in
$\wtI$.

Since we already had
 $$
\bb1-q_1=v_1^*v_1\sim v_1v_1^*=\bb1-p_1\;,
 $$
and since $p_2-q_0\le\bb1-q_0$, we conclude that $p_2-q_0\sim p_0$
for some projection $p_0\le\bb1-p_1$. But then
 $$
p=p_1-w_{12}w_{12}^*+p_0
 $$
is a projection in $A$ and
 $$
p\sim(p_1-w_{12}w_{12}^*)\oplus(p_2-w_{21}w_{21}^*)\;.
 $$
\hfill$\square$
\enddemo

\medskip

Since we will prove later that all extremally rich $C^*-$algebras with
weak cancellation also have $K_1-$surjectivity, the use of $K_1-$surjectivity
in the hypothesis of the next lemma is just a temporary expedient.

\proclaim{3.11. Lemma} Let $A$ be an extremally rich $C^*$--algebra
and $I$ a closed ideal of $A$. Assume that both $I$ and $A/I$ have
weak cancellation and that $eAe/eIe$ has $K_1$--surjectivity for
every projection $e$ in $A$. Then $A$ has weak cancellation.
\endproclaim

\demo{Proof} Let $p$ and $q$ be projections in $A$ which generate the
same closed ideal $J$ in $A$ and have the same image in $K_0(J)$.
Since weak cancellation passes to ideals we may replace $A$ by $J$,
i\.e\. we may assume that $p$ and $q$ are full projections. If
$\pi:A\to A/I$ denotes the quotient map then the conditions above
are also satisfied for $\pi(p)$ and $\pi(q)$ relative to $\pi(A)$ and
$K_0(\pi(A))$. By hypothesis there is therefore an element $u$ in
$pAq$ such that
 $$
\pi(uu^*)=\pi(p)\quad\text{and}\quad
\pi(u^*u)=\pi(q)\;.
 $$
Since $pAq$ is extremally rich by Proposition 3.2, and
$\pi(u)\in\Cal E(\pi(pAq))$, we can apply [{\bf 9}, Theorem 4.1]
and choose $u$ in $\Cal E(pAq)$. (Alternatively use
the argument (i) $\Rightarrow$ (ii) in [{\bf 9}, Theorem 6.1] on
$pAq$.) Let
 $$
p_1=p-uu^*\;,\quad q_1=q-u^*u\;,\quad e_1=u^*u\;.
 $$
Then $p_1$ and $q_1$ belong to $I$ and we can apply Lemma 3.4 to
obtain a full projection $e_2$ in $A$, projections $p_2$ and $q_2$ in
$I$ such that $[p_2]=[q_2]$ in $K_0(I)$, and such that
 $$
p\sim p_2\oplus e_2\;,\quad q= q_2+e_2\;.
 $$

Next, an argument similar to part of the proof of (iii)$\Rightarrow$(i)
in Theorem 3.5 yields
 $$
p\sim p_3\oplus e_3\;,\quad q=q_3 +e_3\;,
 $$
where $e_3$ is a full projection in $A$ and $p_3$, $q_3$ are projections
in $\Cal D(I)$ such that $[p_3]=[q_3]$ in $K_0(\Cal D(I))$, as follows:

\noindent Let $\rho:I\to I/\Cal
D(I)$ denote the quotient map.  From $[\rho (p_2)]=[\rho (q_2)]$ and
$\tsr(\rho (I))=1$ we conclude that $\rho (p_2)\sim \rho (q_2)$.  Then there
is $v$ in $\Cal E(p_2Iq_2)$ such that $\rho(v)$ implements this
equivalence.  Thus
 $$
p_2\sim p_3'\oplus e_3'\;,\quad q_2=q_3'+e_3'\;,
 $$
where $e_3'=v^*v$, $p_3'=p_2-vv^*$, $q_3'=q_2-v^*v$, and
$p_3'\,,q_3'\; \in \Cal D(I)$.  Since $[p_3']-[q_3']$ is in the kernel of
$\iota_*:K_0(\Cal D(I))\to K_0(I)$, there is $\alpha$ in
$K_1(I/\Cal D(I))$ with $\partial_1 \alpha=[p_3']-[q_3']$.  Since $e_2$ is
full in $A$, $e_2Ie_2$ is a full hereditary $C^*-$subalgebra of $I$
and $\rho(e_2Ie_2)$ is full in $I/\Cal
D(I)$.  Thus $\alpha =[\rho(v_0)]$, where $\rho(v_0)$ is unitary in
$\rho(e_2\wtI e_2)$, and $v_0$ may be taken in $\Cal E(e_2\wtI e_2)$.
Hence
 $$
p_3'\oplus e_2\sim p_3\oplus e_3^{''}\;,\quad q_3'+e_2=q_3+e_3^{''}\;,
 $$
where $e_3^{''}=v_0^*v_0$, $p_3\sim p_3'\oplus (e_2-v_0v_0^*)$, $q_3=q_3'+e_2-v_0^*v$, 
and $e_3^{''}$ is full in $A$.  Now let $e_3=e_3'+e_3^{''}$.

Finally we note that $\Cal D(I)=(\bigcup K_j)^=$, where $\{K_j\}$ is an
upward directed family of ideals each of which is generated by finitely
many defect projections $e_3-ww^*$, $w\in \Cal E(e_3\wtI e_3)$.  
Here we are identifying $\wtI$ with $I+\Bbb C\bb1_{\wtA}$.  Since
$K-$theory is compatible with direct limits, and since every projection
in $\Cal D(I)$ is contained in some $K_j$, there are $j_0$ and a finite
collection, $w_1,\dots w_n$, in $\Cal E(e_3\wtI e_3)$ such that
$e_3-w_1w_1^*,\dots,e_3-w_nw_n^*$ generate $K_{j_0}$, $p_3\,, q_3\in K_{j_0}$, 
and $[p_3]=[q_3]$ in $K_0(K_{j_0})$.  Now an argument similar to part of
the proof of (iii)$\Rightarrow$(i) in [{\bf 12}, Lemma 3.8] shows that
there is a projection $e_4$ in $e_3\Cal D(I)e_3$ which generates $K_{j_0}$
(as an ideal).  Since $I$ has weak cancellation, it follows that
$p_3\oplus e_4\sim q_3+e_4$, whence $p\sim p_3\oplus e_3\sim q_3+e_3=q$.
(Since $p_3\oplus e_4\precsim p$, we are not here assuming stable
weak cancellation for $I$.)
\hfill$\square$
\enddemo

\medskip

\example{3.12. Remarks} 
\medskip

(i) It follows that if $A$ is extremally rich and for some $n$,
$\Cal D^n(A)$ is either $0$ or satisfies the hypothesis of
Corollary 3.7, then $A$ has weak cancellation.  Of course, in the latter
case $\Cal D^m(A)=\Cal D^{n+1}(A)$ for $m\ge n+1$.
\medskip

(ii) The $C^*-$algebras called $\Cal E_n$ by Cuntz in [{\bf 16}], for 
$2\le n<\infty$, satisfy the hypotheses of Corollary 3.7 and also 
$\Cal D(\Cal E_n)=\Cal E_n$, but $\Cal E_n$ is not purely infinite, since 
$\Cal E_n$ contains an ideal isomorphic to $\Bbb K$.
\medskip

(iii) If it were known that $A$ extremally rich and $\Cal D(A)=A$ implies
weak cancellation, then it would be easy to deduce that all extremally
rich $C^*-$algebras have weak cancellation.
\medskip

(iv)  It is easy to extend the methods of Cuntz in \cite{\bf 17} to
show that any $C^*-$algebra in which every projection is properly
infinite has weak cancellation, and in fact this is actually a semigroup
result, cf\. [{\bf 21}, Lemma 2.1] and [{\bf 15}, Theorem 4.17].
\medskip

(v)  Since $u^n\in \Cal E(\wtB)$ whenever $u\in \Cal E(\wtB)$, and since
 $$
\bb 1-u^n(u^n)^*\sim (\bb 1-uu^*)\oplus \dots \oplus (\bb 1-uu^*),
 $$
Theorem 3.5(iii) implies that the extremally rich $C^*-$algebra $A$ has
weak cancellation if the set of equivalence classes of projections in
$B$ is ``convex'' in a suitable sense for all $B$ of the
form $pAp$.  Any extra hypotheses
that guarantee this
convexity imply additional positive results.
\endexample
\bigskip

\proclaim{3.13. Corollary} If $A$ is an extremally rich
$C^*$--algebra with weak cancellation, then $\wtA$ has
weak cancellation.
\endproclaim
\hfill$\square$

\bigskip

\proclaim{3.14. Proposition} If $I$ is a closed ideal in an
extremally rich $C^*$--algebra $A$ and if $A$ has weak cancellation, then $A/I$ has
weak cancellation.
 \endproclaim

\demo{Proof} Let $p$ be a projection in $A/I$ and put $B=p(A/I)p$.
Then let $C$ be the inverse image of $B$ in $A$. Since $C$ is a
hereditary $C^*$--subalgebra of $A$ it has weak cancellation and is
extremally rich. (But if the projection $p$ does not lift we may
never find a unital substitute for $C$.) It follows that extreme
points lift from $B$ to $\wtC$, so if $\{u_1,\dots,u_n\}$ is a finite
subset of $\Cal E(B)$ it can be lifted to a finite subset
$\{w_1,\dots,w_n\}$ in $\Cal E(\wtC)$. Since $\wtC$ has weak
cancellation by Corollary 3.13 we can find a projection $q$ in $\wtC$
(actually in $\Cal D(\wtC)=\Cal D(C))$ such that
$q\sim\bigoplus(\bb1-w_iw_i^*)$. It follows that
 $$
\pi(q)\sim\bigoplus(p-u_iu_i^*)\;,
 $$
where $\pi: \wtC \to B$ is the quotient map, 
whence $A/I$ has weak cancellation by Theorem 3.5. \hfill$\square$
\enddemo

\bigskip

\proclaim{3.15. Proposition} If $A$ is an extremally rich
$C^*$--algebra with weak cancellation, then:
\medskip

(i) There is for every projection
$p$ in $\Cal D(A\otimes\Bbb K)$ ($=\Cal D(A)\otimes\Bbb K$) a
projection $q$ in $\Cal D(A)$ such that $p\sim q$.
\medskip

(ii) If $p$ is a projection in $\Cal D(A)\otimes\Bbb K$, then there is an 
infinite set $\{p_n\}$ of mutually orthogonal projections in $\Cal D(A)$ 
such that $p_n\sim p, \forall n$.
\medskip

(iii) If $\Cal D(A)$ is $\sigma-$unital, then $\Cal D(A)$ has a full,
hereditary, stable, $\sigma-$unital $C^*-$subalgebra $B$.
\endproclaim

\demo{Proof} (i) By Corollary 3.13 we may assume that $A$ is unital. Since
$\Cal D(A)$ is generated as an ideal by the set
 $$
\Cal D=\{\bb1-uu^*\mid u\in\Cal E(A)\}
 $$
it follows that $\Cal D(A\otimes\Bbb K)$ is generated by the set
$\Cal D\otimes e_{11}$. There is therefore a finite subset $\{u_i\}$
in $\Cal E(A)$ such that
 $$
p\precsim\bigoplus(\bb1-u_iu_i^*)\;.
 $$
Applying condition (iii) of Theorem 3.5 with $B=A$ we find a
projection $q_0$ in $A$ with $q_0\sim\oplus(\bb1-u_iu_i^*)$. Thus
$p\sim q\le q_0$ and evidently
 $$
q\in\Cal D(A\otimes\Bbb K)\cap A=\Cal D(A)\;.
 $$

(ii) We do a recursive construction.  At step $n$ we construct $n+1$
mutually orthogonal and equivalent projections, $p_1,\dots,p_n,\, q_n$.  The
first step is done by applying part (i) with $2p$ in place of $p$.  At step
$n+1$ we apply part (i) to $B= \her(\bb1-s_n)$, where $s_n=p_1+\dots+p_n$.
Since $B$ is full in $A$, $\Cal D(B) =B\cap \Cal D(A)$, and thus 
$q_n\in \Cal D(B)$.  Hence, we can obtain $p_{n+1}$ and $q_{n+1}$
by applying part (i) with $2q_n$ in place of $p$.

(iii) The hypothesis implies that there is a countable set, $\{p_n\}$, of 
projections in $\Cal D(A)$ such that $\Cal D(A)=\id(\{p_n\})$.  Then the
same technique as in part (ii) produces a countable set, $\{q_m\}$, of
mutually orthogonal projections which consists of infinitely many equivalent
copies of each $p_n$.  Then take  $B=\her(\{q_m\})$.
\hfill$\square$
\enddemo
\smallskip

\example{Remark}  Of course, we cannot require $B$ to be a corner. It is 
possible that $A$ is unital and $\Cal D(A)=A$.
\endexample 
\medskip

\proclaim{3.16. Corollary} If $A$ is an extremally rich $C^*-$algebra with
weak cancellation, then $A$ also has stable weak cancellation.
\endproclaim

\demo{Proof} We apply Theorem 3.5 to prove that $\Cal D(A)\otimes\Bbb K$
has weak cancellation.  This is immediate since each subalgebra of
the form $p(\Cal D(A)\otimes\Bbb K)p$ is isomorphic to an algebra
$q\Cal D(A)q$.  Then Lemma 3.11 implies $A\otimes \Bbb K$ has weak
cancellation, since $\tsr((A\otimes \Bbb K)/(\Cal D(A)\otimes\Bbb K))=1$.
\hfill$\square$
\enddemo

\bigskip

The next corollary fulfills a promise made in [{\bf 12}, Remark 3.6].

\proclaim{3.17. Corollary}  If $A$ is an extremally rich $C^*-$algebra
with weak cancellation and if $p$ is a projection in $A$ such that
$\Cal D(pAp)=pAp$ (equivalently $\Cal D(I)=I$ where $I=id(p)$), then
$p$ is properly infinite.  
In particular, the hypotheses of Corollary 3.7 imply that every defect
projection is properly infinite.
\endproclaim

\demo{Proof}  By applying the Proposition to $pAp$, we immediately conclude
that $2p \precsim p$.  \hfill$\square$
\enddemo

\bigskip

\proclaim{3.18. Proposition} If $A$ is a unital $C^*$--algebra which
is extremally rich and has weak cancellation there is for each $n$
and every $u$ in $\Cal E(\Bbb M_n(A))$ a $v$ in $\Cal E(A)$ such that
 $$
\bb1-vv^*\sim\bb1_n-uu^*\quad\text{and}\quad
\bb1-v^*v\sim\bb1_n-u^*u\;.
 $$
\endproclaim

\demo{Proof} By Proposition 3.15 we can find projections $p_1$ and
$q_1$ in $\Cal D(A)$ such that
 $$
2(\bb1_n-uu^*)\sim p_1\quad\text{and}\quad
2(\bb1-uu^*)\sim q_1\;.
 $$
Thus for subprojections $p\le p_1$ and $q\le q_1$ we have
 $$
\bb1_n-uu^*\sim p\precsim\bb1-p\quad\text{and}\quad
\bb1_n-u^*u\sim q\precsim\bb1-q\;.
 $$
But then both $\bb1-p$ and $\bb1-q$ generate $A$ as an ideal, and
evidently $[\bb1-p]=[\bb1-q]$ in $K_0(A)$ (since $[p]=[q]$). By weak
cancellation $\bb1-p=vv^*$ and $\bb1-q=v^*v$ for some partial
isometry $v$ in $A$. But since $\bb1_n-uu^*$ and $\bb1_n-u^*u$ are
centrally orthogonal in $\Bbb M_n(A)$ we see that $p$ and $q$ are
centrally orthogonal in $A$, whence $v\in\Cal E(A)$. \hfill$\square$
\enddemo

\bigskip

\proclaim{3.19. Corollary} With $A$ as above, if $\Bbb M_n(A)$
contains a proper isometry so does $A$. (Finiteness implies 
stable finiteness.)
\endproclaim
\bigskip

\example{3.20. Remarks} 
\medskip

(i) One may ask whether a unital and extremally rich
$C^*$--algebra $A$ contains a proper isometry if $\Bbb M_n(A)$ does
for some $n$, and whether this happens if $\Cal D(A)=A$. In the 
presence of weak cancellation the answer is yes in both cases by 
Corollaries 3.17 and 3.19. Some condition 
beyond $\tsr(A)>1$ is necessary, though,
even when $A$ has weak cancellation. See the discussion of infiniteness
conditions in [{\bf 12}] and note that the extended Toeplitz algebra 
has weak cancellation by Lemma 3.11. 
\medskip

(ii)  Propositions 3.17 and 3.18 also imply that if $A$ satisfies the
hypotheses of Corollary 3.7, then so does $\Bbb M_n(A)$.

\endexample
\bigskip

\subhead{4. $K_1$--Surjectivity}\endsubhead

\medskip

In the next three lemmas we shall be concerned with a closed ideal $I$
in a unital $C^*$--algebra $A$ and $\wtI$ will denote $I+\Bbb
C\bb1$. (The possibility $I=A$ is not excluded.) In $\Bbb M_2(A)$ we
consider the unital $C^*$--subalgebra $B$ consisting of matrices of
the form
 $$
\pmatrix
a & x_{12}\\
x_{21} & \lambda\bb1+x_{24}\endpmatrix\;,\quad
\lambda\in\Bbb C\;,\;x_{ij}\in I\;,\;a\in A\;.
 $$
Note that the subset of $B$ determined by $\lambda =0$ is an 
ideal of $B$ which is Rieffel-Morita equivalent to $A$.  (If
$I=A$ this ideal is the whole of $B$).  Thus
every ideal $J$ of $A$ gives rise to an ideal of $B$ (which is
determined by $a\in J$, $x_{ij}\in I\cap J$, and $\lambda =0$).  
We shall commit a
slight abuse of notation and denote both ideals by the same symbol.
We shall denote by $B_{00}^{-1}$ the connected component containing
$\bb1$ of the group of invertible elements in 
$\Bbb M_2(\wtI )\cap B$.  

\proclaim{4.1. Lemma} Assume $\wtI$ has weak cancellation and
$v\,,w\in \wtI$.
If $v\in\Cal E(\wtI)$ and
$v^*v+w^*w=\bb1$, and if furthermore $\left(\smallmatrix w\\
                                     v\endsmallmatrix\right)$
is the second column of a left invertible element of $B$,
then $\bb1-ww^*\sim vv^*$ in $\wtI$.
\endproclaim

\demo{Proof} Since $v\in\Cal E(\wtI)$ we know that $w$ is a partial
isometry $(w^*w=\bb1-v^*v)$ and
 $$
[\bb1-ww^*]=[\bb1]-[w^*w]=[v^*v]=[vv^*]
 $$
in $K_0(\wtI)$. Thus the lemma follows by weak cancellation if we can
show that both projections generate $\wtI$ as an ideal.

In fact it is
equivalent to show that they generate $A$ as an ideal.  The condition
for a projection $p$ in $\wtI$ to generate $\wtI$ as an ideal is twofold:

\, (i) Either $p\notin I$ or $I=A$.

(ii) The algebra $pIp$ generates $I$ as an ideal.

\noindent But if $p$ generates $A$ as an ideal, then (i) is obviously true, and also
$pAp$ is a full hereditary $C^*-$subalgebra of $A$.  And from this we
easily deduce that $pIp$ $(=pAp\cap I)$ is full in $I$.

Now the fullness of $vv^*$ 
follows from Lemma 3.3. 
Let $J$ denote the closed ideal of $A$ generated by $\bb1-ww^*$. If
$J\ne A$ we pass to $A/J$ without changing notation. The conditions
on $v$ and $w$ are unchanged, but now $ww^*=\bb1$. Since $v\in\Cal
E(A)$ we have $\bb1-vv^*$ centrally orthogonal to $w^*w$
($=\bb1-v^*v$), and this now forces $\bb1-vv^*=0$. Thus both $v$ and
$w$ are co--isometries. But this contradicts the last part of the
hypothesis.
\hfill$\square$
\enddemo

\bigskip

\proclaim{4.2. Lemma} Suppose that $b$ is a left invertible element
of $B$ and that $\wtI$ is extremally rich with weak cancellation.
There is then an element $b_0$ in $B_{00}^{-1}$ such that $b_0b=\pmatrix
a & 0\\ 0 & \bb1\endpmatrix$ for some left invertible element $a$ in $A$.
\endproclaim

\demo{Proof} Write $b=\pmatrix \star & y\\ \star & x\endpmatrix$. By
assumption $b^*b$ is invertible so $y^*y+x^*x\in(\wtI)_+^{-1}$. If
$y^*=u_1|y^*|$ is the polar decomposition of $y^*$ and $f$ is a
function in $C_0((0\,,\Vert y\Vert ])$ such that $|f(t)t-t|$ is small for all 
$t$ in $[0\,,\Vert y\Vert ]$ then with $y_1=f(|y^*|)|y^*|u_1^*$ we still have
$y_1^*y_1+x^*x\in(\wtI)_+^{-1}$. By [{\bf 9}, Theorem 3.3] there is an
element $v_1$ in $\Cal E(\wtI)$ such that $v_1^*y_1+x$ is
quasi--invertible. Thus $v_1^*y_1+x=ev$ for some $v$ in $\Cal
E(\wtI)$ and $e$ in $(\wtI)_+^{-1}$, cf\. [{\bf 9}, Theorem 1.1].
Multiplying $b$ from the left by the element
 $$
\pmatrix \bb1 & 0\\
0 & e^{-1}\endpmatrix\;
\pmatrix \bb1 & 0\\
v_1^*f(|y^*|) & \bb1\endpmatrix
 $$
 in $B_{00}^{-1}$ gives a matrix $\pmatrix \star & y\\ \star &
v\endpmatrix$. Left multiplication by $\pmatrix \bb1 & -yv^*\\ 0
& \bb1\endpmatrix$, also in $B_{00}^{-1}$, gives the matrix $\pmatrix
\star & z\\ \star & v\endpmatrix$, where $z=y(\bb1-v^*v)$. Since we
still have $z^*z+v^*v$ invertible and $v^*v$ is a projection this
implies that $z^*z\ge\varepsilon (\bb1-v^*v)$ for some $\varepsilon
>0$. If therefore $z=w|z|$ is the polar decomposition of $z$ then
$w\in\wtI$ (actually $w\in I$) and we consider $wh(|z|)$ in $I$ for
some $h$ in $C_0((0\,,\Vert z\Vert ])$ such that $h(|z|)|z|=\bb1-v^*v$. Left
multiplication with
 $$
\pmatrix
\bb1-ww^*+h(|z^*|) & 0\\
0 & \bb1
\endpmatrix\;,
 $$
which belongs to $B_{00}^{-1}$, transforms the matrix $\pmatrix \star &
z\\ \star & v\endpmatrix$ into $\pmatrix \star & w\\ \star &
v\endpmatrix$, where now $w^*w=\bb1-v^*v$.

The elements $v,w$ now satisfy the conditions in Lemma 4.1 since
$\pmatrix \star & w\\ \star & v\endpmatrix$ is left invertible. Thus
$v^*v=w_1^*w_1$ and $\bb1-ww^*=w_1w_1^*$ for some partial isometry
$w_1$ in $\wtI$. Since $v\in\Cal E(\wtI)$ it has the form
$v=\lambda\bb1+c$ for some $c$ in $I$ and $\lambda$ in $\Bbb C$ with
$|\lambda|=1$. After left multiplication with $\pmatrix \bb1 & 0\\ 0
& \overline\lambda\endpmatrix$ we can assume that $\lambda=1$ without
any other changes. Since $w^*w+w_1^*w_1=\bb1$ and
$ww^*+w_1w_1^*=\bb1$ the element $u=w+w_1$ is unitary in $\wtI$. For
$0\le t\le 1$ define
 $$
b_t=\pmatrix
\bb1 & -u\\
0 & \bb1\endpmatrix\;
\pmatrix
\bb1 & 0\\
t(\bb1-v)u^* & \bb1\endpmatrix\;
\pmatrix \bb1 & w_1v^*\\
0 & \bb1\endpmatrix\;.
 $$
Since $\bb1-v\in I$ it follows by routine calculations that $b_t\in
B_{00}^{-1}$ for all $t$. 
However,
 $$
b_1\pmatrix
\star & w\\
\star & v\endpmatrix\;=\;
\pmatrix
a & 0\\
d & \bb1\endpmatrix
 $$
for some elements $a$ in $A$ and $d$ in $I$. Since the original
element $b$ was left invertible it now follows that $a$ must be left
invertible (in $A$); so if $a_1a=\bb1$ we perform the final left
multiplication by the element
 $$
\pmatrix
\bb1 & 0\\
-da_1 & \bb1\endpmatrix\quad\text{in}\;\;B_{00}^{-1}
 $$
to obtain the desired solution $\pmatrix a & 0\\ 0 &
\bb1\endpmatrix$. \hfill$\square$
\enddemo

\bigskip

\proclaim{4.3. Lemma} Suppose that $b$ is a quasi--invertible
element in $B$ and that $\wtI$ is extremally rich with weak
cancellation. We can then find $b_1$ and $b_2$ in $B_{00}^{-1}$ such
that $b_1bb_2=\pmatrix a & 0\\ 0 & \bb1\endpmatrix$ for some
quasi--invertible element $a$ in $A$.
\endproclaim

\demo{Proof} Without changing notation, we replace $A$ with 
$A\bigoplus \Bbb C$, $I$ with $I\bigoplus 0$, and $b$ with
$b\oplus \lambda\bb 1$.  This avoids an annoying but essentially trivial
complication that would occur if $I\ne A$ but $(I+J_2)/J_2=A/J_2$.

By [{\bf 9}, Theorem 1.1] there are orthogonal closed
ideals $J_1$ and $J_2$ of $B$ with corresponding quotient morphisms
$\pi_1$ and $\pi_2$ such that $\pi_1(b)$ is left invertible and
$\pi_2(b)$ is right invertible. Using Lemma 4.2 with $\pi_2(B)$,
$\pi_2(I)$ and $\pi_2(b^*)$ in place of $B$, $I$ and $b$, which is
legitimate by Proposition 3.14, we find an element $\pi_2(b_2)$ in
$\pi_2(B)_{00}^{-1}$ such that
 $$
\pi_2(b)\pi_2(b_2)=
\pmatrix
\pi_2(a_2) & 0\\
0 & \bb1\endpmatrix\;.
 $$
Since invertible elements in the connected component of the identity
are always liftable we may assume that $b_2\in B_{00}^{-1}$, and we can
write
 $$
bb_2=\pmatrix a_2 & x_{12}\\ x_{21} & \bb1+x_{22}\endpmatrix\;,
\quad \text{with}\quad x_{ij} \quad \text{in} \quad I_2 = I \cap J_2\;.
 $$

Define the $C^*$--algebra
 $$
B_2=\pmatrix
A & I_2\\
I_2 & \wtI_2\endpmatrix\;\subset B
 $$
and note that the restriction of the morphism $\pi_1$ to $B_2$ is an
isomorphism except possibly at the $(1,1)$--corner since 
$J_1\cap J_2=\{0\}$.
As $\pi_1(I_2)$\;$(=I_2$) has weak cancellation, we
can apply Lemma 4.2 with $\pi_1(B_2)$, $\pi_1(I_2)$ and $\pi_1(bb_2)$
in place of $B$, $I$ and $b$; and find an element $\pi_1(b_1)$ in
$\pi_1(B_2)_{00}^{-1}$ such that
 $$
\pi_1(b_1)\pi_1(bb_2)=
\pmatrix
\pi_1(a_1) & 0\\
0 & \bb1\endpmatrix\;.
 $$
Again we may assume that $b_1\in(B_2)_{00}^{-1}$, but since $\pi_1|B_2$
is an isomorphism except at the $(1,1)$--corner this implies that
 $$
b_1bb_2=\pmatrix
a & 0\\
0 & \bb1\endpmatrix
 $$
for some $a$ in $A$. Necessarily then $a\in A_q^{-1}$, as desired.
\hfill$\square$
\enddemo

\bigskip

\proclaim{4.4. Theorem} Every extremally rich $C^*$--algebra with
weak cancellation has $K_1$--surjectivity.
\endproclaim

\demo{Proof} By Corollary 3.13 we may assume that the $C^*$--algebra
$A$ is unital. We can therefore 
use Lemma 4.2 for (two--sided) invertible elements and $I=A$ in
every morphism
 $$
\Bbb M_{2^n}(A)^{-1} \to \Bbb M_{2^{n+1}}(A)^{-1}\;.
 $$
\hfill$\square$
\enddemo

\medskip

\proclaim{4.5. Theorem} In the category of extremally rich
$C^*$--algebras the subcategory of algebras that also have 
weak cancellation is stable under the formation of quotients,
hereditary $C^*$--subalgebras (in particular ideals), matrix
tensoring, Rieffel-Morita equivalence, arbitrary extensions, 
and inductive limits.  Also if the extremally rich $C^*-$algebra
$A$ has a composition series of ideals, $\{I_{\alpha}\;|\;0\le \alpha \le \lambda\}$,
such that $I_0=0$, $I_{\lambda}=A$, and $I_{\alpha +1}/I_{\alpha}$ has weak
cancellation for each $\alpha < \lambda$, then $A$ has weak cancellation.
\endproclaim

\demo{Proof} The wording of the result reflects the fact that the
category of extremally rich $C^*$--algebras is not itself stable
under (arbitrary) extensions, cf\. [{\bf 9}, Theorem 6.1], and 
stable only under extreme point preserving inductive limits, cf\. 
[{\bf 9}, Proposition 5.2].

To verify the claims: Quotients follow from Proposition 3.14;
hereditary is trivial; 
extensions follow from Lemma 3.11 in
conjunction with Theorem 4.4; and inductive limits respect ideals,
$K$--theory and equivalence.  Corollary 3.16 now implies that weak
cancellation is stable under both matrix tensoring and tensoring with
$\Bbb K$, hence under stable isomorphism.  By [{\bf 8}] Rieffel-Morita
equivalence coincides with stable isomorphism when both algebras are
$\sigma -$unital.  The general case can be reduced to the separable case
by the same technique as in [{\bf 9}, Lemma 5.6 and Theorem 5.7].
For the last sentence we prove by transfinite induction that each $I_{\alpha}$
has weak cancellation. If $\alpha$ is a limit ordinal, then $I_{\alpha}$ is
a direct limit, and otherwise it is an extension.
\hfill$\square$
\enddemo

\medskip

\bigskip

\definition{4.6. The Extremal $K-$set, $K_e$}  In [{\bf 11}] 
we defined two analogues of $K_1$ which use extremal partial
isometries in place of unitaries, or equivalently, quasi-invertibles
in place of invertibles.  One of these, denoted $\Cal E_{\infty}(A)$,
takes two extremals, each in some matrix algebra over $\wtA$, to be
equivalent if $u\oplus \bb 1_k$ is homotopic to $v\oplus \bb 1_l$
in $\Cal E(\Bbb M_n(\wtA ))$ for some (large) $n$ and suitable $k\,,l$.
The equivalence relation for $K_e(A)$ is coarser and is given in
[{\bf 11}, Definitions 3.6].  For example, the defect ideals,
$I= \text{id}(\bb 1 -uu^*)$ and $J= \text{id}(\bb 1 -u^*u)$, are invariants
of the $K_e-$class of $u$, as are also the classes of $\bb 1 -uu^*$
and $\bb 1 -u^*u$ in $K_0(I)$ and $K_0(J)$, respectively.  But even
the Murray-von Neumann equivalence classes of $\bb 1 -uu^*$ and $\bb 1 -u^*u$
are invariants of the $\Cal E_{\infty}-$class of $u$.  If this were
the only difference between $K_e$ and $\Cal E_{\infty}$,  then obviously weak
cancellation would imply $K_e=\Cal E_{\infty}$.  Although the difference
is more extensive, we shall prove in the next section that
$K_e(A)=\Cal E_{\infty}(A)$ when $A$ is extremally rich with weak
cancellation.  Neither $K_e$ nor $\Cal E_{\infty}$ is a group, but both
contain $K_1$ and the group $K_1$ acts on both.  
The next result is that extremal richness with weak
cancellation implies $K_e-$surjectivity.  The same proof shows
``$\Cal E_{\infty}-$surjectivity'', a property which is formally stronger,
but equivalent in this situation. 
\enddefinition
\bigskip

\proclaim{4.7. Theorem} Every extremally rich $C^*-$algebra with
weak cancellation has $K_e-$surjectivity.
\endproclaim

\demo{Proof}  The proof is similar to that of Theorem 4.4, except that
we use Lemma 4.3 instead of 4.2.  \hfill$\square$
\enddemo

\bigskip

\subhead{5. Good Index Theory}\endsubhead

\medskip

\proclaim{5.1. Theorem} Every extremally rich $C^*-$algebra with 
weak cancellation has good index theory.
\endproclaim

\demo{Proof} We are given a unital $C^*-$algebra $A$, an ideal $I$ which
is extremally rich with weak cancellation, the quotient map
$\pi :A\to A/I$, a unitary $\overline u$ in $A/I$, and $\alpha$ in $K_1(A)$ such
that $\pi_*(\alpha)=[\overline u]$ in $K_1(A/I)$.  We wish to find a unitary
$u$ in $A$ such that $\pi (u) =\overline u$.  Let $\pi_n$ denote the
natural map from $\Bbb M_n(A)$ to $\Bbb M_n(A/I)$.  We may choose $n$,
a power of $2$,
so that there is a unitary $v$ in $\Bbb M_n(A)$ which belongs to the
class $\alpha$ such that $\pi_n(v)$ is homotopic to 
$\overline u \oplus \bb1_{n-1}$ in $\Cal U (\Bbb M_n(A/I))$.  
Then $(\overline u \oplus \bb1_{n-1})(\pi_n(v))^{-1}$ can be lifted to
$w$ in $\Cal U_0(\Bbb M_n(A))$.  We replace $v$ with $wv$, without changing
notation, and thus achieve that $\pi_n(v)=\overline u \oplus \bb1_{n-1}$.

Thus $v$ belongs to the algebra called $B$ in connection with Lemma 4.2,
with $\Bbb M_{n/2}(A)$ in place of $A$, and that Lemma provides $b_0$ in
$B_{00}^{-1}$ such that $b_0b$ has the form $v' \oplus \bb 1_{n/2}$. 
The $(1,1)-$corner of $b_0$ is congruent to a scalar modulo $I$, and
clearly we may assume this scalar is $1$.  Thus $v'$ satisfies the same
properties as $v$, relative to $n/2$, except for the inconsequential fact
that it is only invertible instead of unitary.  We may remedy this, if
desired, with a polar decomposition.  Continuing in this way, we attain
our goal.  \hfill$\square$
\enddemo
\bigskip

\example{Remark}  The proof actually provides $u$ in the given class
$\alpha$.  Thus we could have dispensed with Theorem 4.4 and proved
$K_1-$surjectivity and good index theory simultaneously.  We hope the
reader will forgive this and a few other minor inefficiencies
in the organization of the paper.
\endexample
\bigskip

\proclaim{5.2. Proposition}  Assume $u$ and $v$ are extremal partial 
isometries in matrix algebras over the unital $C^*-$algebra $A$ which lie in
the same $K_e-$class.  If the defect ideals are extremally rich with
weak cancellation, then $u$ and $v$ lie in the same $\Cal E_{\infty}-$class.
\endproclaim

\demo{Proof} Let $I$ and $J$ be the left and right defect ideals, which
are the same for $u$ and $v$.  (Here we regard the defect ideals as ideals
of $A$, using the identification of ideals of $A$ with ideals of
$\Bbb M_n(A)$.)  By replacing $A$ with $\Bbb M_n(A)$ for suitable $n$,
$u$ with $u\oplus \bb 1_{n-k}$, $v$ with $v\oplus \bb 1_{n-l}$, and
changing notation, we may assume $u$ and $v$ are in $A$.  Let
$\pi:A\to A/(I+J)$, $\rho:A\to A/J$, and $\lambda:A\to A/I$ be the
quotient maps.

Since $\pi (w)=\pi (v)\pi (u)^{-1}$ is a unitary whose class in $K_1(A/(I+J))$
is $0$, it follows from Theorems 4.4 and 5.1 that we may take $w$ to be
a unitary whose class in $K_1(A)$ is $0$.  Thus $\pi(v)=\pi(u')$,
where $u' =wu$.  We now construct unitaries $w_1\in \bb 1+I$ and
$w_2\in \bb 1+J$ such that the classes of $w_1$ and $w_2$ in
$K_1(I)$ and $K_1(J)$, respectively, are trivial, $\rho(v)=\rho(w_1u')$,
and $\lambda(v)=\lambda(u'w_2)$.  Once this is done, we have that
$v=w_1u'w_2=w_1wuw_2$, since $I\cap J=0$.  Since all of the $w$'s 
are trivial in $K_1(A)$, it follows that $u$ and $v$ are equivalent in
$\Cal E_{\infty}(A)$.

To construct $w_1$, we may replace $A$ with $A/J$, since $\rho$ is an
isomorphism on $I$.  Then $u'$ and $v$ are isometries which agree
modulo $I$.  The defect projections $p=\bb 1-u'u'{}^*$ and 
$q=\bb 1-vv^*$ each generate the ideal $I$ and have the same class
in $K_0(I)$.  Thus there is $x$ such that $x^*x=p$ and $xx^*=q$.
Then $w'=x+vu'{}^*$ is a unitary in $\bb 1+I$, and $v=w'u'$.  Now since
$K_1(pIp)=K_1(I)$ and $pIp$ has $K_1-$surjectivity, there is $y$ in
$\Cal U(pIp)$ which induces the same class as $w'$ in $K_1(I)$.  So we
can take  $w=w'(\bb 1-p+y^*)$.  The construction of $w_2$ is similar.
\hfill$\square$
\enddemo

\bigskip

\proclaim{5.3. Corollary} If $A$ is extremally rich with weak cancellation,
then $K_e(A)=\Cal E_{\infty}(A)$.  \hfill$\square$
\endproclaim
\bigskip
\definition{5.4.Extremal Analogues of Good Index Theory}

Four different
properties are listed below.  In all cases $K$
is an ideal in a unital $C^*-$algebra $A$ and
$\pi_n :\Bbb M_n(A) \to \Bbb M_n(A/K)$ are the quotient maps.
\roster
\item If $\overline u\in \Cal E(A/K)$, $v\in \Cal E(\Bbb M_n(A))$, and 
$[\pi_n(v)]_{K_e}=[\overline u]_{K_e}$, then there is $u\in \Cal E(A)$
such that $\pi(u)=\overline u$ and $[u]_{K_e}=[v]_{K_e}$.
\item If $\overline u\in \Cal E(A/K)$, $v\in \Cal E(\Bbb M_n(A))$, and 
$[\pi_n(v)]_{\Cal E_{\infty}}=[\overline u]_{\Cal E_{\infty}}$, then there is $u\in \Cal E(A)$
such that $\pi(u)=\overline u$ and $[u]_{\Cal E_{\infty}}=[v]_{\Cal E_{\infty}}$.
\item If $\overline u\in \Cal E(A/K)$, $v\in \Cal E(\Bbb M_n(A))$, and 
$[\pi_n(v)]_{K_e}=[\overline u]_{K_e}$, then there is $u\in \Cal E(A)$
such that $\pi(u)=\overline u$.
\item If $\overline u\in \Cal E(A/K)$, $v\in \Cal E(\Bbb M_n(A))$, and 
$[\pi_n(v)]_{\Cal E_{\infty}}=[\overline u]_{\Cal E_{\infty}}$, then there is $u
\in \Cal E(A)$
such that $\pi(u)=\overline u$.
\endroster

Obviously $(1)$ implies $(3)$ and $(2)$ implies $(4)$, but we cannot
assert, for example, that $(1)$ (for all choices of $A$)
is equivalent to $(3)$ plus $K_e-$surjectivity (for $K$),
because we have no exact sequence controlling the lack of injectivity
of $K_e(\pi)$.  (Exception:  If $\tsr(K)=1$, then [{\bf 11}, Theorem 5.4]
fills this gap and implies the equivalence of $(1)$ and $(2)$ with
$(3)$ and $(4)$, respectively.)  Also $(3)$ implies $(4)$, since it
derives the same conclusion from a weaker hypothesis.  But there is no
obvious comparison between $(1)$ and $(2)$, since both the hypothesis and
conclusion of $(1)$ are weaker.  We shall prove that $(2)$ is true
whenever $K$ is extremally rich with weak cancellation.  We don't know
whether $(1)$ is true under the same hypothesis or, for example, whenever
$K$ is extremally rich of real rank zero.  But by applying also Proposition 5.2,
we see that $(1)$ is true whenever $K$ and the defect ideals of
$\overline u$ are extremally rich with weak cancellation.  Also, an
easy pullback argument shows that $(1)$ is true if
$K$ is extremally rich with weak cancellation and
$(I+J)\cap K=0$, where $I$ and $J$ are the defect ideals of $v$.
Finally, it can be shown that if $(3)$ or $(4)$ holds for $A=M(K)$,
the multiplier algebra, then it holds for all choices of $A$.  The proof
is based on Busby's analysis of extensions, [{\bf 14}].  The technical
issue that arose in [{\bf 13}, \S 4], namely that the map
$\tau:A/K\to M(K)/K$ may not be extreme-point-preserving, causes no
dificulty here.  In fact, the existence of $v$ implies that $\tau(\overline u)$
is in $\Cal E(M(K)/K)$.

\enddefinition
\bigskip
\proclaim{5.5. Theorem} Let $K$ be a closed ideal in a unital
$C^*$--algebra $A$ and assume that $K$ is extremally rich with 
weak cancellation. If $\overline u\in \Cal E(A/K)$, 
$v\in \Cal E(\Bbb M_n(A))$, and
$[\pi_n(v)]_{\Cal E_{\infty}}=[\overline u]_{\Cal E_{\infty}}$, where 
$\pi_n :\Bbb M_n(A) \to \Bbb M_n(A/K)$ is the quotient map, 
then there is $u\in \Cal E(A)$
such that $\pi_1(u)=\overline u$ and 
$[u]_{\Cal E_{\infty}}=[v]_{\Cal E_{\infty}}$.
\endproclaim

\demo{Proof} The proof is almost identical to that of Theorem 5.1, except
that we use Lemma 4.3 instead of 4.2.  We may assume that $n$ is a power
of $2$ and that $\pi_n(v)$ is homotopic to $\overline u\oplus \bb1_{n-1}$
in $\Cal E (\Bbb M_n(A/K))$.  By [{\bf 11}, Corollary 2.3] there are
$\pi_n(w_1)$ and $\pi_n(w_2)$ in $\Cal U_0 (\Bbb M_n(A/K))$
such that $\overline u\oplus \bb1_{n-1}=\pi_n(w_1)\pi_n(v)\pi_n(w_2)$. 
We may assume $w_1, w_2\in \Cal U_0 (\Bbb M_n(A))$.  Then without
changing notation, we replace $v$ with $w_1vw_2$ to achieve
$\pi_n(v)=\overline u\oplus \bb1_{n-1}$.  Then coninue as in 5.1.
\hfill$\square$
\enddemo

\bigskip

\proclaim{5.6. Corollary} 
Let $K$ be a closed ideal in a unital
$C^*$--algebra $A$ and assume that $K$ is extremally rich with
weak cancellation. If $\overline u\in \Cal E(A/K)$,
$v\in \Cal E(\Bbb M_n(A))$, and
$[\pi_n(v)]_{K_e}=[\overline u]_{K_e}$, where
$\pi_n :\Bbb M_n(A) \to \Bbb M_n(A/K)$ is the quotient map, and if the defect ideals of $\overline u$ are extremally rich with weak cancellation, then there is $u$ in $\Cal E(A)$ such that $\pi_1(u)=\overline u$ and $[u]_{K_e}=[v]_{K_e}$. 
\endproclaim

\demo{Proof} By Proposition 5.2 the hypotheses of the Corollary imply
those of the Theorem.
\hfill$\square$
\enddemo

\bigskip
\proclaim{5.7. Corollary} 
Let $K$ be a closed ideal in a unital
$C^*$--algebra $A$ and assume that $K$ is extremally rich with
weak cancellation. If $\overline u\in \Cal U(A/K)$,
$v\in \Cal E(\Bbb M_n(A))$, the defect ideals of $v$ are in $K$, and
$[\pi_n(v)]_{K_1}=[\overline u]_{K_1}$, where
$\pi_n :\Bbb M_n(A) \to \Bbb M_n(A/K)$ is the quotient map,
then there is $u\in \Cal E(A)$
such that $\pi_1(u)=\overline u$ and
$[u]_{K_e}=[v]_{K_e}$.  \hfill$\square$
\endproclaim
\bigskip
\proclaim{5.7$'$. Corollary} Let $K$ be a closed ideal in a unital
$C^*$--algebra $A$ and assume that $K$ is extremally rich with
weak cancellation. If $\overline u\in \Cal U(A/K)$,
and if $[\overline u]_{K_1}$ is in the image of $K_e(A)$, then $u$ can
be lifted to $\Cal E(A)$.       \hfill$\square$
\endproclaim

\bigskip
\proclaim{5.8. Proposition} Let $K$ be a closed ideal in a unital
$C^*$--algebra $A$ and assume that $K$ is extremally rich with
weak cancellation. If $\overline u\in \Cal E(A/K)$,
$v\in \Cal E(\Bbb M_n(A))$, and
$[\pi_n(v)]_{K_e}=[\overline u]_{K_e}$, where
$\pi_n :\Bbb M_n(A) \to \Bbb M_n(A/K)$ is the quotient map,
and if $(I+J)\cap K=0$, where $I$ and $J$ are the defect ideals of
$v$,
then there is $u\in \Cal E(A)$
such that $\pi_1(u)=\overline u$ and
$[u]_{K_e}=[v]_{K_e}$.
\endproclaim

\demo{Proof}  Consider the pullback diagram
$$
\CD
A @>{\rho}>> A/(I+J)\\
@V{\pi}VV  @VV{\tau}V \\
A/K @>>{\sigma}> A/(I+J+K)
\endCD
$$
where the maps are the obvious ones and $\pi=\pi_1$.  
Then $\sigma(\overline u)$ is a
unitary whose $K_1-$class is lifted by $[\rho_n(v)]\in K_1(A/(I+J))$.
By Theorems 5.1 and 4.4 there is $w$ in $\Cal U(A/(I+J))$ such that
$\tau(w)=\sigma(\overline u)$ and $[w]=[\rho_n(v)]$ in 
$K_1(A/(I+J))$.  Thus there is $u$ in $A$ such that $\pi(u)=\overline u$
and $\rho(u)=w$, whence $u\in \Cal E(A)$.  Now the defect ideals of $u$
are contained in $I+J$ and map under $\pi$ to the defect ideals of
$\overline u$, the defect ideals of $\overline u$ are the same as those
of $\pi_n(v)$, and $\pi_{|I+J}$ is an isomorphism.  Thus $u$ also has
defect ideals $I$ and $J$.  Then it follows from
[{\bf 11}, Theorem 4.5] that $[u]_{K_e}=[v]_{K_e}$.  \hfill{$\square$}

\enddemo
\bigskip
\example{5.9. Remarks}

(i) The relations of the results in this section to classical index theory
become clearer if reformulated in an equivalent way.  Thus for
good index theory we would start with $x$ in $A$ such that $\pi(x)$ is
invertible in $A/K$ (such an $x$ is called a Fredholm element relative to
$K$) and seek a $K-$perturbation of $x$ which is invertible.  Similarly,
for the extremal analogues of good index theory we would start with
$x$ such that $\pi(x)\in (A/K)_q^{-1}$ (such an $x$ was called quasi-Fredholm
in [{\bf 11}, Definitions 6.3]) and seek a $K-$perturbation in
$A_q^{-1}$.  Quasi-Fredholm elements are meant to be analogous to
classical semi-Fredholm operators (but we have used the name
quasi-invertible, not semi-invertible).  Also, instead of hypothesizing 
$v$ in $\Cal E(\Bbb M_n(A))$,  
we would hypothesize a class $\alpha$ in $K_e(A)$.  Furthermore, in
[{\bf 11}, Definitions 6.3] we also defined an index space, $\text{Ind}_e(K)$,
which is the orbit space of $K_e(A/K)$ under the image of $K_1(A)$.  The
existence of $\alpha$ could be reformulated in terms of the index in this
sense of $[\overline u]_{K_e}$ (in 5.4(1) or (3)).  Now every element 
$\overline \alpha$ of
$\text{Ind}_e(K)$ has built into it a pair $(I,J)$ of defect ideals
and a class $\beta$ in $K_0(D)$, where $D=\pi^{-1}(I+J)$, and $\beta$ is
obtained from the boundary map just as in the Fredholm case.  Moreover,
for given $(I,J)$, $\overline \alpha$ is determined by $\beta$.  
(However, it is awkward to describe which classes $\beta$ arise in this
way.)  Since $\beta$ lives in $K_0(D)$ instead of $K_0(K)$, it seems 
reasonable that we should use hypotheses on the defect ideals as well as
on $K$ to prove 5.4(1) or (3).

(ii) Corollaries 5.7 and 5.7$'$ should be compared to a result of
G. Nagy, [{\bf 29}, Theorem 2].  This implies {\it a fortiori}
\medskip

$\quad\qquad$ If $K$ has general stable rank (gsr) at most $2$, and if $\overline u$ is an element

\noindent (N)$\;\;\qquad$        of $\Cal U(A/K)$ such that $\partial_1 ([\overline u]_{K_1}) \le 0$, then $\overline u$ can be lifted to an iso-

$\qquad\quad$ metry in $A$.
\medskip

\noindent Now the hypothesis $\text{gsr}(K) \le 2$ is not comparable with
our hypothesis on $K$, but it is implied by $\tsr(K) =1$ or even
$\text{csr}(K)\le 2$. Aside from this difference, (N) is intermediate 
in strength between 5.7 and 5.7$'$.  It is fairly routine to deduce from
$\partial_1 ([\overline u]_{K_1}) \le 0$ the existence of an isometry $v$
in some $\Bbb M_n(A)$ such that $[v]_{K_e}$ lifts $[\overline u]_{K_1}$.
Thus 5.7 gives the conclusion of (N) and also allows us to control the
$K_e-$class of the lift if $[v]_{K_e}$ is given, whereas 5.7$'$ states only
that $\overline u$ can be lifted to $\Cal E(A)$ and doesn't require that
the lift be an isometry.  It is also interesting that even though our
hypothesis on $K$ doesn't imply Nagy's, nevertheless Nagy's proof will
work with our hypothesis.
\medskip

(iii) If $A/K$ is extremally rich with weak cancellation, then the hypothesis
on defect ideals in Corollary 5.6 is automatically satisfied.
Also, by the last remark in 5.4, if the corona
algebra $C(K)=M(K)/K$ is extremally rich with weak cancellation, then
5.4(3) is true for all $A$ with no hypotheses other than those on $K$.
If the corona algebra hypothesis is satisfied also for all ideals of $K$,
then we even get 5.4(1).  To see this we first use Proposition 5.8, 
applied to $A/(K\cap(I+J))$, to reduce to the case $K\subset I+J$.  Then the
argument based on [{\bf 14}] applies not just to give a lift $u$ but
to show that $u$ has the same defect ideals as $v$.  But then 
[{\bf 11}, Theorem 4.5] implies that $[u]_{K_e}=[v]_{K_e}$.

\endexample
\bigskip

\example{5.10. Example} For ease of notation set $\Bbb B=\Bbb B(H)$
for some infinite dimensional separable Hilbert space $H$ and choose a
projection $p$ in $\Bbb B$ such that both spaces $p(H)$ and 
$(\bb1-p)(H)$ are infinite-dimensional. Let $A$ be the
$C^*$--subalgebra of $\Bbb B\otimes c$ consisting of convergent
sequences $x=(x_n)$ such that
 $$
\lim (\bb1-p)x_np=\lim px_n(\bb1-p)=0\;.
 $$
This algebra was considered in [{\bf 9}, Examples 1.3 \& 5.3] to give
an example of a $C^*$--algebra which is the inductive limit of
extremally rich $C^*$--algebras (actually von Neumann algebras)
without itself being extremally rich.

Let $I=\Bbb B\otimes c_0$, which is clearly a closed ideal in $A$,
and consider the (split) extension
 $$
0 @>>> I @>>> A @>>> \Bbb B\oplus\Bbb B @>>> 0\;.
 $$
In this piquant situation all $K$--groups vanish; but the extremal
$K$--sets do not, and they control the quasi--Fredholm elements in
$A$ since $I$ and $A/I$ are extremally rich
with weak cancellation. The Fredholm theory is trivial in
this example: Every invertible element in $A/I$ lifts to a
invertible element in $A$ -- as it must by Theorem 5.1.

Writing $\Bbb Z^e=\Bbb Z\cup\{\pm\infty\}$ we find that $K_e(I)$ is
the set of sequences in $\Bbb Z^e$ that are eventually zero, whereas
$K_e(A/I)=(\Bbb Z^e)^2$. An element $x=(x_n)$ belongs to $A_q^{-1}$  
if and only if every $x_n$ is either left or right invertible and there is an
$\varepsilon >0$ such that for all $n$, $|x_n|$ (and $|x_n^*|$)
has a gap $]0,\varepsilon [$ in its spectra, cf\. [{\bf 9}, Theorem
1.1]. Since $(x_n)$ converges to a block diagonal operator in $\Bbb
B^2$ we can describe $K_e(A)$ as eventually constant sequences
$(\alpha_n)$ in $\Bbb Z^e$ together with an element
$(\alpha_\infty^1,\alpha_\infty^2)$ in the first or third quadrant of
$(\Bbb Z^e)^2$ such that
$\alpha_\infty^1+\alpha_\infty^2=\lim\alpha_n$. Elements $(\alpha_n)$
and $(\beta_n)$ in $K_e(A)$ are composable (c\.f\. [{\bf 11}, \S 2.6]) if and only if $\alpha_n$ and
$\beta_n$ have the same sign for $1\le n\le\infty$, in which case
$(\alpha_n)+(\beta_n)=(\alpha_n+\beta_n)$.

Thus the image of $K_e(A)$ in $K_e(A/I)$ is $(\Bbb Z_+^e)^2\cup(\Bbb
Z_-^e)^2$, and a quasi-Fredholm element can be perturbed to an
element of $A_q^{-1}$ if and only if its index is in the union of the first and
third quadrants.  Equivalently, an element of $\Cal E(A/I)$ can be
lifted to $\Cal E(A)$ if and only if both components of its $K_e-$class
have the same sign.  

In view of [{\bf 11}, \S 7] it is also interesting to consider
$K=\Bbb K\otimes c_0$ and $K_1= \{x\in A\;|\;x_n\in \Bbb K, \forall n\}$.
Then $K=\overline{Soc(A)}$, and $K_1/K=\overline{Soc(A/K)}$.  Here the
ordinary $K-$groups do not all vanish, and it can be seen that
$K_e(A/K)=K_e(A/K_1)=\{((\alpha_n),\alpha_{\infty}^1,\alpha_{\infty}^2):\;
\alpha_n, \alpha_{\infty}^1,\alpha_{\infty}^2\in \Bbb Z^e\,,\,
(\alpha_{\infty}^1,\alpha_{\infty}^2)\ne(\infty,-\infty),(-\infty,\infty)\,,\,
\alpha_n=\alpha_{\infty}^1 + \alpha_{\infty}^2\text{ , eventually}\}$.
\endexample
\bigskip

\proclaim {5.11. Proposition} Suppose that $A$ is an extremally rich
$C^*$--algebra with weak cancellation such that the natural map 
 $$
\Cal U(M(A))/\Cal U_0(M(A))\quad @>>> \quad K_1(M(A))
 $$
is injective. Then also the following map is injective:
 $$
\Cal U(C(A))/\Cal U_0(C(A)) \quad @>>> \quad K_1(C(A)) \;.
 $$

\endproclaim

\demo{Proof} 
Consider the 
commutative diagram
 $$
\CD
\Cal U(M(A)) @>>{[\quad]}> K_1(M(A))\\
@VV{\pi}V                  @VV{\pi_1}V\\
\Cal U(C(A)) @>{[\quad]}>> K_1(C(A))\;.
\endCD
 $$

If $u$ is in $\Cal U(C(A))$ such that $[u]=0$ then by Theorems 5.1 and 4.4
we have $u=\pi(v)$ for some $v$ in $ \Cal U(M(A))$ with $[v]=0$.
By assumption this means that $v \in \Cal U_0(M(A))$, whence 
$u \in \Cal U_0(C(A))$. \hfill$\square$
\enddemo

\bigskip

The next result does not mention extremal richness, but it illustrates 
an application of good index theory.

\proclaim {5.12. Corollary} If $A$ is $\sigma$--unital and stable, then we 
have a short exact sequence of groups
 $$
0 @>>> \Cal U_0(C(A)) @>>> \Cal U(C(A)) @>>> K_0(A) @>>> 0 \;.
 $$
\endproclaim

\demo{Proof} By the Cuntz-Higson-Mingo result, [{\bf 18}], [{\bf 28}] or 
[{\bf 39}, Theorem 16.8], the group
$\Cal U(M(A))$ is connected (even contractible), so exactness at 
$\Cal U(C(A))$ follows from the proof of Proposition 5.11 and the 
previously mentioned results from [{\bf 29}], [{\bf 30}], [{\bf 35}], and [{\bf 38}]. Exactness at 
$K_0(A)\quad (\cong K_1(C(A)))$ is well known. It follows from the fact 
that $\bb1$ is equivalent to $\bb1_n$ in $\Bbb M_n(C(A))$. \hfill $\square$

\enddemo

\bigskip

\example {5.13. Remark} There are many other cases for which it is known 
that $\Cal U(M(A))$ is connected, and, of course, the above arguments can apply
to these as well. For example, Theorem 2.4 of [{\bf 20}] states that 
$\Cal U(M(A))$ is connected when $A$ is a separable, matroid 
$C^*$--algebra, and Elliott's proof works equally well for arbitrary 
$\sigma$--unital AF--algebras. Lin proves in [{\bf 26}, Lemma 3.3] that 
$\Cal U(M(A))$ is connected for some $C^*$--algebras of real rank zero 
and stable rank one. The hypothesis that tsr$(A)=1$ is used only to 
ensure (strong) cancellation, and it seems obvious that in some cases 
weak cancellation would suffice.
\endexample

\bigskip

\subhead{6. $K_1-$injectivity and $K_0-$surjectivity}\endsubhead

\medskip
The main goal of this section is to prove that extremal richness plus
weak cancellation 
implies $K_1-$injectivity.  This is accomplished in
Theorem 6.7, the main step being Lemma 6.5, which already
includes all the extremally rich $C^*-$algebras which are
purely properly infinite. 
Our proof is partly modeled on Cuntz's $K_1-$injectivity proof in [{\bf 17}],
but we need some additional ideas, in particular the introduction of
$K_0-$surjectivity.  We also use a technique similar to one used by
Zhang in [{\bf 41}].

\definition{6.1. The map $\partial_0:K_0(A/I)\to K_1(I)$}
Since Bott periodicity identifies $K_0(A/I)$ with $K_1(\text{S}(A/I))$,
where $\text{S}$ denotes suspension, we may consider $\partial_0$ to be
defined on this latter group. We need to know the form of $\partial_0 \beta$
in the special case $\beta=[u]$, $u\in \Cal U(\widetilde{\text{S}(A/I)})$. 
Thus, $u$ is given by a continuous function
$f:[0,1]\to \Cal U(\widetilde{A/I})$ such that $f(0)=f(1)=\bb 1$.  Then $f$ can
be lifted to $g:[0,1]\to \Cal U(\wtA)$ such that $g(0)=\bb 1$ and
$g(1)\in (\bb1+I)\cap \Cal U(\wtI)$, and $\partial_0 \beta=[g(1)]$.  It
is important to note that $g(1)$ is null-homotopic in $\Cal U(\wtA)$.
\enddefinition
\medskip

\definition{6.2. Definitions} We say that $A$ has {\it (strong) 
$K_0-$surjectivity} if the group $K_0(A)$ is generated by
$\{[p]\;|\;p\text{ is a projection in }A\}$.
Thus Zhang's result
in [{\bf 40}] shows that $C^*-$algebras of real rank zero have
strong $K_0-$surjectivity.  
In [{\bf 17}] Cuntz showed that purely infinite simple $C^*-$algebras 
satisfy a still stronger property, which, however, is too strong for our 
purposes below. 
We say that $A$
has {\it weak} $K_0-${\it surjectivity} if $\text{S}A$ has $K_1-$surjectivity.  Then strong $K_0-$surjectivity implies
weak $K_0-$surjectivity because the function defined by
$f(t)=\text{exp}(2\pi itp)$ is a unitary in $(\text{S}A)\sptilde$ which
corresponds to $p$ under Bott periodicity, for each projection $p$ in $A$.
Since Rieffel showed in [{\bf 35}] that $\text{csr}(A)\le 2$ implies
$K_1-$surjectivity for $A$ and also that $\tsr(A)\le 1$ implies
$\text{csr}(\text{S}A)\le 2$, we see that all $C^*-$algebras of stable rank
one have weak $K_0-$surjectivity.
\enddefinition
\medskip

\proclaim{6.3. Proposition} If $A$ is extremally rich with weak cancellation,
then $\Cal D(A)$ has strong $K_0-$surjectivity.
\endproclaim

\demo{Proof} Since $\Cal D(A)$ is generated as an ideal by projections,
$K_0(\Cal D(A))$ is generated by $\{[p]\;|\;p\text{ is a projection in }\Cal D(A)\otimes \Bbb K\}$.  But by Proposition 3.15, each such $p$ is equivalent
to a projection in $\Cal D(A)$.		\hfill$\square$
\enddemo

\proclaim{6.4. Lemma} If $A$ is a $C^*-$algebra
and $B$ is a $\sigma-$unital hereditary $C^*-$subalgebra such that
$B^{\perp}$ contains an infinite set $\{p_n\}$ of mutually orthogonal
and equivalent projections which are full (in $A$),
then there is a full hereditary 
$C^*-$subalgebra $B'\supset B$ such that 
$B'\cong B'{}'\otimes \Bbb K$, where $B'{}'$ is unital.  In
particular, $B'$ has an approximate identity,
$(e_n)$, consisting of full projections, such that  for each $m$ and $n$, 
$me_n\precsim \bb1_{\wtA}-e_n$. 
\endproclaim

\demo{Proof} 
Let $C=\her(p_1,p_2,\dots)$ and $B'=\text{her}(C\cup B)$.  
Then since $C\cong (p_1Ap_1)\otimes \Bbb K$, and $B'$ is $\sigma-$unital,
[{\bf 7}, Theorem 4.23] implies that $B'\cong C$.  
(The Kasparov Stabilization Theorem, [{\bf 23}, Theorem 2] could also
be used for the last part of the proof.  As explained on page 963 of
[{\bf 7}], [{\bf 6}, Theorem 3.1] (the main ingredient of 
[{\bf 7}, Theorem 4.23]) and [{\bf 23}, Theorem 2] are essentially
equivalent.)
\hfill$\square$
\enddemo
\bigskip

\proclaim{6.5. Lemma} If $A$ is an extremally rich $C^*$--algebra 
with weak cancellation, then $\Cal D(A)$ has 
$K_1-$injectivity.
\endproclaim

\demo{Proof} Let $D=\Cal D(A)$ and let $\wtD$ be the forced unitization.
Assume there is $u$ in $(\bb1+D)\cap \Cal U(\wtD)$ whose 
$K_1-$class is trivial  but $u$ is not null-homotopic 
in $\Cal U(\wtD)$.  We claim then that there is an ideal $J$ of $D$ which
is maximal with respect to the property that $u+J$ fails to be null-homotopic
in $\Cal U(\wtD /J)$.  To prove this by Zorn's Lemma, we may assume a
totally ordered collection $\{J_i\}$ of ideals such that, with
$J=(\bigcup J_i)^=$, $u+J$ is null-homotopic and prove that for some $i$,
$u+J_i$ is null-homotopic in $\Cal U(\wtD /J_i)$.
Since $\Cal U_0(\wtD/J)$ is the image of $\Cal U_0(\wtD)$,
$u$ is homotopic in $\Cal U(\wtD)$ to some $v \in (\bb1+J)\cap \Cal U(\wtJ)$.
Writing $v=\bb1+x$, $x\in J$, we see that for some $i$,
$\Vert x+J_i\Vert <1$, whence $v+J_i$ is null-homotopic.
Now by [{\bf 9}, Theorem 6.1], extremal partial isometries lift from
$\wtA/J$ to $\wtA$, and thus $\Cal D(A/J)=D/J$.  Therefore we may replace
$A$ with $A/J$ and $D$ with $D/J$, without changing notation, and seek
to obtain a contradiction.

Next choose a continuous function $f:\Bbb T\to [0,\infty)$, where 
$\Bbb T$ is the unit circle, such that 
$\{z\;|\;f(z)\ne 0\}=\{e^{i\theta}\;|\;2\pi/3<\theta <4\pi/3\}$.
Let $B_1=\her(f(u))$, and note that $B_1\subset D$, since $f(1)=0$.
If $B_1=0$, then the spectrum of $u$ omits $-1$, a contradiction.

Case (i):  If $\tsr(B_1)=1$, let $I=\id(B_1)$.
By construction,
$u+I$ is null-homotopic.  Hence we see as above that $u$ is homotopic
to some $v \in  \Cal U(\wtI)$.  If $\alpha$ is the class of $v$ in $K_1(I)$, 
then $\alpha$ maps to
$0$ in $K_1(D)$.  It follows that $\alpha =\partial_0 \beta$ for some
$\beta$ in $K_0(D/I)$.  Now Proposition 6.3 implies that $D/I$ has 
$K_0-$surjectivity (since $\Cal D(A/I)=D/I$).  Thus 6.1 applies, and we
see that $\alpha$ is represented by a unitary w in $\wtI$ which is 
null-homotopic in $\Cal U(\wtD)$.  But now $w^*v$ is a unitary in $\wtI$ 
whose $K_1-$class vanishes and $\tsr(I)=1$ (for example by
[{\bf 9}, Corollary 5.8]).  Thus Rieffel's $K_1-$ injectivity result
in [{\bf 35}] now applies to give our contradiction.

Case (ii):  If $\tsr(B_1)>1$, then let $p$ be a non-zero defect projection
for $B_1$ and note that $\bb1-p$ is full (by Lemma 3.3, for example).
Let $I=\id(p)$, and let $\pi:\wtD \to \wtD /I$ be the quotient map.
Since $\text{Re}(pup)\le -\tfrac12 p$, $pup$ is invertible in $pDp$
and is homotopic to $p$ within $(pDp)^{-1}$.
It follows that $u$ is homotopic to $p+u_1$, where $u_1$ is a unitary element
of
$(\bb 1-p)\wtD (\bb 1-p)$ such that $\pi(u_1)=\pi(u)$.  To see this,
first homotop $u$ within ${\wtD}^{-1}$ to
 $$
(\bb 1-(\bb1-p)u(pup)^{-1})u(\bb1-(pup)^{-1}u(\bb1-p)),
 $$
which has the form $pup+u'$, $u'$ invertible in $(\bb 1-p)\wtD (\bb 1-p)$.

By construction $\pi(u_1)$ is null-homotopic.  Since $\pi(\wtD)=\pi((\bb 1-p)\wtD (\bb 1-p))$,
every element of $\Cal U_0(\wtD/I)$ lifts to
$\Cal U_0((\bb 1-p)\wtD (\bb 1-p))$.  It follows that $u_1$ is homotopic
(within $\Cal U((\bb 1-p)\wtD (\bb 1-p))$) to a unitary $u_2$ in
$(\bb 1-p)\wtI (\bb 1-p)$.

If $\alpha$ is the class of $u_2$ in $K_1(I)$, then $\alpha$ maps to
$0$ in $K_1(D)$.  It follows that $\alpha =\partial_0 \beta$ for some
$\beta$ in $K_0(D/I)$.  As above we use the $K_0-$surjectivity of $D/I$ 
and 6.1 to get a special form for $\partial_0 \beta$, but now we use
$(\bb 1-p)\wtD (\bb 1-p)$ in place of $\wtD$.  So we find a 
unitary $v$ in $(\bb 1-p)\wtI (\bb 1-p)$ such that $v$ is null-homotopic
in $\Cal U((\bb1-p)\wtD(\bb1-p))$ and $[v]=\alpha$ in $K_1(I)$.  Then if
$u_3=u_2v^*$, we see that $p+u_3$ is homotopic to $u$ in $\Cal U(\wtD)$
and $[p+u_3]=0$ in $K_1(I)$.  
From now on all the action takes place in $\wtI$, and we will make
no further use of the assumption that $u$ has a null-homotopic image in
any non-trivial quotient.

Next let $\rho:\wtI \to \wtI/\Cal D(I)$ be the quotient map.  By [{\bf 35}],
which we apply
within $\rho((\bb1-p)\wtI(\bb1-p))$, $\rho(u_3)$ is null-homotopic.
(Here we are using the fullness of $\bb1-p$, which implies that
$[u_3]=0$ in $K_1((\bb1-p)I(\bb1-p))$.)  Since $\rho((\bb1-p)I(\bb1-p))$ has weak
$K_0-$surjectivity, we may use the above argument to homotop $u_3$ to an
element $u_4$ in $\Cal U((\bb1-p) \widetilde{\Cal D(I)}(\bb1-p))$ such that 
$[u_4]=0$ in $K_1(\Cal D(I))$.

Then we write $\Cal D(I)=(\bigcup I_j)^=$, where each $I_j$ is the ideal
generated by finitely many defect projections, and $\{I_j\}$ is directed
upward.  Then $u_4$ is homotopic to 
$u_5$ in some $\Cal U((\bb1-p)\widetilde{I_j}(\bb1-p))$, and because of the compatibility
of $K_1$ with direct limits, we may assume $[u_5]=0$ in $K_1(I_j)$.
Clearly we may also assume $u_5\in\bb1-p+(\bb1-p)I_j(\bb1-p)$.
Since $p$ is full in $I$, $\Cal D(pIp)=pIp\cap \Cal D(I)$, a full hereditary
$C^*-$subalgebra of $\Cal D(I)$.  Hence every projection in $\Cal D(I)$ is
equivalent to a projection in $\Cal D(pIp)\otimes \Bbb K$.
So Proposition 3.15(ii) can be applied to $pIp$ to find an infinite set, 
$\{q_n\}$, of mutually orthogonal and mutually equivalent projections in
$pIp$, which are full in $I_j$.  (So that actually $q_n\in pI_jp$.)

Finally we apply Lemma 6.4 with $I_j$ in place of $A$  and 
$B=\text{her}(u_5-\bb 1+ p)$.  If $B'$ and $\{e_n\}$ are as in the
Lemma, let $t_n=\bb 1 +e_n(u_5-\bb 1+ p)e_n$.  Since $t_n\to p+u_5$,
for large $n$ there is a unitary $w_n$ in $e_nI_je_n$ such that
$\bb 1 -e_n +w_n$ is homotopic to $p+u_5$ in $\Cal U(\wtI_j)$.  Since 
$e_n$ is full in $I_j$, $K_1(e_nI_je_n)$ is naturally isomorphic to
$K_1(I_j)$, and hence $[w_n]=0$ in $K_1(e_nI_je_n)$.  It follows that for
some $m$,
$w_n \oplus me_n$ is null-homotopic in $\Cal U(\Bbb M_{m+1}(e_nI_je_n))$.
But the conclusion of Lemma 6.4 states that $me_n\precsim \bb 1 -e_n$.
Thus $\bb 1 -e_n+w_n$ is null-homotopic in $\Cal U(\wtI_j)$, and $u$ is
null-homotopic in $\Cal U(\wtD)$.
\hfill$\square$
\enddemo
\bigskip

\proclaim{6.6. Proposition} Let $I$ be an ideal of a $C^*-$algebra
$A$.  Then:
\medskip
{\rm (i)} If $I$ and $A/I$ have $K_1-$injectivity and $A/I$ has
weak $K_0-$surjectivity, then $A$ has $K_1-$injectivity.
\medskip

{\rm (ii)} If $I$ and $A/I$ have weak $K_0-$surjectivity and $I$ has
$K_1-$injectivity, then $A$ has weak $K_0-$surjectivity.
\endproclaim
\medskip

\demo{Proof} We may assume $A$ unital.  Let $\pi:A\to A/I$ be the
quotient map.

The argument for part (i) already occurred several times
in the proof of Lemma 6.5.  If $[u]_{K_1}=0,\,u \in \Cal U(A)$, then
$\pi(u)$ is null-homotopic, whence there is $v$ in $\Cal U(\wtI)$ which
is homotopic to $u$.  If $\alpha=[v]_{K_1(I)}$, then $\alpha=\partial_0\beta$, 
and 6.1 applies to give $w$ in $\Cal U(\wtI)$ which is null-homotopic in 
$\Cal U(A)$ and which represents the class $\alpha$.  Then the
$K_1-$injectivity of $I$ is applied to $w^*v$.

For part (ii) we are given $\alpha$ in $K_0(A)$, which is identified with
$K_1(\text{S}A)$.  Then $K_0(\pi)(\alpha)$ is represented by $\overline u$ 
in $\Cal U(\widetilde{\text{S}(A/I)})$, and $\overline u$ is represented by 
a function $f:[0,1]\to \Cal U(A/I)$ such that $f(0)=f(1)=\bb1$.  Lift f to
$g:[0,1]\to \Cal U(A)$ such that $g(0)=\bb1$.  Then, since
$\partial_0(K_0(\pi)(\alpha))=0$, we see from 6.1 that $[g(1)]_{K_1(I)}=0$
and so $g(1)$ is null-homotopic in $\Cal U(\wtI)$. 
Thus there is $h:[0,1]\to \Cal U(\wtI)$ such that $h(0)=\bb1$ and
$h(1)=g(1)$, so that $h^*g$ gives an element of $\Cal U(\widetilde{\text{S}A})$ 
which represents a class $\beta$ in $K_0(A)$.  Finally, $\alpha-\beta$ is
in the image of $K_0(I)$ (which is the kernel of $K_0(\pi)$) and is
therefore represented by a unitary.
\hfill$\square$
\enddemo
\bigskip

\proclaim{6.7. Theorem}  If $A$ is an extremally rich $C^*-$algebra with
weak cancellation, then $A$ has $K_1-$injectivity and weak
$K_0-$surjectivity.
\endproclaim

\demo{Proof} By Lemma 6.5 and Proposition 6.3, $\Cal D(A)$ has both 
properties, and since $\tsr(A/\Cal D(A))=1$, $A/\Cal D(A)$ also has both
properties.  Now apply Proposition 6.6.		\hfill$\square$
\enddemo
\bigskip

\proclaim{6.8. Theorem}  If $A$ is an extremally rich $C^*-$algebra
with weak cancellation, then $A$ also has
$K_e-$injectivity.
\endproclaim

\demo{Proof} Let $u$ and $v$ be elements of $\Cal E(\wtA)$ which lie in
the same $K_e-$class.  Exactly as in the proof of Proposition 5.2, we
show that $v=w_1wuw_2$, where all the $w$'s are unitaries in $\wtA$
whose $K_1-$classes vanish.  Thus, by 
Theorem 6.7, all the $w$'s are in $\Cal U_0(\wtA)$,
and hence $u$ is homotopic to $v$ in $\Cal E(\wtA)$.
\hfill$\square$
\enddemo
\bigskip
The next result should be regarded as an example.

\proclaim{6.9. Corollary} If $A$ is a purely infinite simple $C^*-$algebra,
then any two proper isometries in $\wtA$ are homotopic within the
set of isometries.
\endproclaim

\demo{Proof} By Cuntz's results in [{\bf 17}] $A$ has weak
cancellation (and $K_1-$injectivity), and by [{\bf 11}, Corollary 4.7] and its proof, all proper
isometries lie in the same $K_e-$class.		\hfill$\square$
\enddemo
\bigskip

The next result is a summary theorem which includes the results we
have proved about three classes of extremally rich $C^*-$algebras,
except that it omits the results related to
the extremal analogues of good index theory.

\proclaim{6.10. Theorem} If $A$ is an extremally rich $C^*-$algebra, then,
under any one of the hypotheses listed below,
$A$ has weak cancellation, $K_1-$bijectivity (i\.e\.,
$K_1(A)=\Cal U(\wtA)/\Cal U_0(\wtA)$), good index theory, 
$K_e-$bijectivity (i\.e\., $K_e(A)=\Cal E(\wtA)/\text{homotopy}$),
and weak $K_0-$surjectivity.
Moreover, $\Cal D(A)$ has strong $K_0-$surjectivity.
\medskip

{\rm (i)}  $A$ has real rank zero.
\smallskip

{\rm (ii)} If $I$ is the left defect ideal of an
element of $\Cal E(\wtA)$, then $\Cal D(I)=I$.  In particular, this applies
if $A$ is purely infinite or if every defect projection is properly infinite.
\smallskip

{\rm (iii)} $A$ is isometrically rich.
\endproclaim

\demo{Proof}  Combine 
Theorems 3.10, 4.4, 4.7, 5.1, 6.7, and 6.8,
Corollaries 3.6 and 3.7, and Proposition
6.3.
\hfill$\square$
\enddemo
\bigskip
\subhead{7. Additional Results and Remarks}\endsubhead

\bigskip

\example{7.1. The Type I Case and the Almost Hausdorff Case}
\endexample
\bigskip

\proclaim{7.1.1.  Lemma} If $A$ is an extremally rich $C^*-$algebra whose
primitive ideal space is Hausdorff, then $A$ has weak
cancellation.
\endproclaim

\demo{Proof} We verify condition (iii) in Theorem 3.5.  Thus let $B=pAp$
for a projection $p$ in $A$, and consider $u\in \Cal E(B)$ with left and
right defect ideals $I$ and $J$.  Since $I^{\vee}$ is a compact-open
subset of $B^{\vee}$, which is Hausdorff, we see that $I^{\vee}$ is
closed.  Thus $B=I\oplus I^{\perp}$.  Let $u=v\oplus w$.  From
$J\subset  I^{\perp}$ it follows that $v$ is an isometry and $w$ a co-isometry.
Hence if $u'=v\oplus \bb1_{I^{\perp}}$, then $u'$ is an isometry with
the same left defect projection as $u$.  Now if $u_1, u_2,\dots, u_n$
are in $\Cal E(B)$, then $s=\prod_1^nu_i'$ is an isometry in $B$ and
 $$
p-ss^*\sim \bigoplus_{i=1}^n(p-u_iu_i^*)
 $$
\hfill$\square$
\enddemo
\bigskip

\proclaim{7.1.2 Proposition}  If $A$ is an extremally rich $C^*-$algebra
whose primitive ideal space is almost Hausdorff, then $A$ has weak
cancellation.
\endproclaim

\demo{Proof} The hypothesis implies that $A$ has a composition series
of ideals, $\{I_{\alpha}\;|\;0\le \alpha \le \lambda\}$, such that
$I_0=0$, $I_{\lambda}=A$, and $(I_{\alpha +1}/I_{\alpha})^{\vee}$ is
Hausdorff  for each $\alpha < \lambda$.  Thus the result follows
from the Lemma and the last sentence of Theorem 4.5.    \hfill$\square$
\enddemo
\bigskip

\proclaim{7.1.3. Corollary}  If $A$ is an extremally rich $C^*-$algebra
which is of type I, then $A$ has weak cancellation.
\hfill$\square$
\endproclaim
\bigskip

\example{7.2. More on Extensions and $K_1-$injectivity} 
\endexample
\bigskip

In this subsection we state with at most
minimal indications of proof some results which are relevant mainly to
non-extremally rich $C^*-$algebras.  Our original plan called
for some of these results to be included in a separate paper to be written by
the first named author.

\proclaim{7.2.1. Theorem} Assume $I$ is an ideal of a $C^*-$algebra $A$,
$I$ and $A/I$ have weak cancellation, $I$ has real rank zero, and
$pAp/pIp$ has $K_1-$surjectivity for each projection $p$ in $A$.  Then
$A$ has weak cancellation.
\endproclaim
\medskip

\proclaim{7.2.2.  Theorem}  Assume $I$ is an ideal of a $C^*-$algebra $A$,
$A/I$ has weak cancellation, $\wtB$ has stable weak cancellation for
every hereditary $C^*-$subalgebra of $I$, and 
$pAp/pIp$ has $K_1-$surjectivity for each projection $p$ in $A$. Then
$A$ has weak cancellation.
\endproclaim

\medskip

Theorem 7.2.1 is closely related to [{\bf 2}, Theorem 7.5], which has a 
similar conclusion when $A$ has real rank zero.
The proof of 7.2.1 is somewhat similar to that of Lemma 3.11, one
difference being that a different method is used for lifting partial
isometries.  The proof of 7.2.2 is also somewhat similar.  Here a
key difference is that the boundary map, 
$\partial_1:K_1(A/I)\to K_0(I)$, is dealt with by the method applicable
to general $C^*-$algebras--i\.e\., the unitary in
$A/I$ need not be liftable to a partial isometry in $A$.

The proof of the next theorem is somewhat similar to parts of the 
proof of Lemma 6.5 and is easier on the whole.

\proclaim{7.2.3. Theorem} Every purely properly infinite $C^*-$algebra
has $K_1-$injectivity.
\endproclaim
\medskip

Just as we found it necessary to link weak cancellation with 
$K_1-$surjectivity and $K_1-$injectivity with $K_0-$surjectivity
to facilitate several proofs, so the following more obvious linkage 
can be useful.

\proclaim{7.2.4. Proposition} Let $I$ be an ideal of a $C^*-$algebra
$A$.  Then:
\medskip
{\rm (i)} If If $I$ and $A/I$ have $K_1-$surjectivity and $I$ has
good index theory, then $A$ has $K_1-$surjectivity.
\medskip

{\rm (ii)} If $I$ and $A/I$ have good index theory and $A/I$ has
$K_1-$surjectivity, then $A$ has good index theory.
\endproclaim
\bigskip

Propositions 6.6 and 7.2.4, and Theorem 7.2.2, combined with easy direct limit arguments, can be
used to derive results for $C^*-$algebras that have composition series
with well-behaved quotients.  Here is one example.
In [{\bf 13}] the authors said that
$A$ has {\it generalized stable rank one} if $A$ has a composition series of
ideals$\{I_{\alpha}\;|\;0\le \alpha \le \lambda\}$
,
such that $I_0=0$, $I_{\lambda}=A$, and $\tsr(I_{\alpha +1}/I_{\alpha})=1$ 
for each $\alpha < \lambda$.  Since [{\bf 13}, Proposition 5.2] 
implies that every type I extremally rich $C^*-$algebra has generalized
stable rank one, 
the following result includes Corollary 7.1.3.

\proclaim{7.2.5. Proposition} Every $C^*-$algebra of generalized stable rank one
has stable weak cancellation, $K_1-$bijectivity, good index theory, and
weak $K_0-$surjectivity.
\endproclaim
\bigskip

\example{7.3. Concluding Remarks}
{}
\endexample
\bigskip

\example{7.3.1. Question}  Does every extremally rich $C^*-$algebra have
weak cancellation?
\medskip

(i) We are {\bf NOT CONJECTURING} either answer to this question, but 
our results suggest that negative examples will not be easy to come by.
\medskip

(ii)  It follows from Theorem 4.5 that an extremally rich $C^*-$algebra
$A$ has a largest ideal $I$ with weak cancellation, and that $B=A/I$ has no
non-zero ideals with weak cancellation.  Then also $B$ has no hereditary
$C^*-$subalgebras with weak cancellation, and hence for every hereditary
$C^*-$subalgebra $C$, there is $u\in\Cal E(\wtC)$ such that neither
defect projection, $p$ or $q$, vanishes.  Then one can apply the
same reasoning with $C$ replaced by $pBp$ or $qBq$ and continue
indefinitely.  
So a negative answer to 7.3.1 implies an example with a rich ideal
structure.  (But see point (iv) below.)
\medskip

(iii)  Several concepts of infiniteness for extremally rich $C^*-$algebras
were discussed in [{\bf 12}].  
The most infinite case is the purely properly infinite
case, and Corollary 3.7 has a weaker hypothesis than this.
We have usually thought of the stable rank one case as the most finite 
(despite that fact that stable rank one algebras can be purely infinite
(c\.f\. [{\bf 37}])).  Also, by
[{\bf 12}, Theorem 3.9], if the extremally rich $C^*-$algebra $A$ is
not purely properly infinite, then there are two ideals,
$J\subsetneqq I$ such that $\tsr(I/J)=1$.  Thus to verify weak
cancellation for $A$, it is sufficient to verify it for $J$ and $A/I$.
But such arguments won't prove anything general, and moreover (see  the
next point) it is illusory to think that they even suggest a positive
answer to 7.3.1.
\medskip

(iv)  In [{\bf 12}, 2.11 and 2.12] we constructed an
extremally rich $C^*-$algebra that exhibits the rich ideal structure
discussed in point (ii).  This algebra, called $B_{II}$, was constructed
with a somewhat different purpose in mind.  It satisfies some
concept of infiniteness, namely it has no non-zero ideal of stable rank one--
by [{\bf 12}, Theorem 2.8] this concept has some superficial resemblence
to Cuntz's definition of purely infinite--but it is stably finite.
The ideal structure of $B_{II}$ makes it impervious to attack (for
the purpose of proving weak cancellation) via extension theory or composition
series.  And the intermediate type of infiniteness that $B_{II}$ has
is also useless for this purpose.  The sequence $(\Cal D^n(B_{II}))$ is strictly decreasing
with intersection $0$.

Nevertheless, it is easy to prove that $B_{II}$ has weak cancellation,
since it is a direct limit of more tractable algebras.  Possibly
extremally rich $C^*-$algebras exist with similar ideal structure that
are not direct limits of more tractable algebras.  But it might be
difficult to construct one explicitly and prove that it is extremally
rich.  Another strategy to construct a counterexample, or to gain
insight from a failed attempt at a counterexample, might be to consider
$A=\varinjlim A_n$, where the $A_n$'s are not extremally rich but
$A$ is (c\.f\. [{\bf 19}]).  Of course, it would be necessary for
$A$ to have a rich ideal structure, unlike the situation in [{\bf 19}].
\medskip

(v) The upshot of all of the above is that there is insufficient evidence
to justify either conjecture for Question 7.3.1.
\endexample
\medskip

\example{7.3.2}
The proof of Theorem 3.10 really shows the following:  
If $A$ is an extremally rich $C^*-$algebra such that for each
projection $p$ in $A$, every quotient of $pAp$ has $K_1-$injectivity
(or even the weaker property that $\Cal U/\Cal U_0$ is commutative),
then $A$ has weak cancellation.  Thus Question 7.3.1 is equivalent to the
question whether every extremally rich $C^*-$algebra has $K_1-$injectivity. 
Why did we devote most of our effort to weak cancellation rather than
$K_1-$injectivity?  Part of the reason is that 
Rieffel's proof of $K_1-$injectivity in [{\bf 35}] is based on his result that
$\text{csr}(\text{S}A)\le \tsr(A) +1$, and we never saw how to 
use the extremal richness of $A$ in a direct way to prove 
anything about $\text{S}A$.
This also explains why we approached weak $K_0-$surjectivity
only in an indirect way.
On the other hand, the idea in Proposition 6.3 to prove 
strong $K_0-$surjectivity for $\Cal D(A)$, not $A$, is clearly correct. 

Of the three parts of Theorem 6.10, case (iii) is the widest in scope,
and it is this case which most justifies the belief that extremal richness
is a useful hypothesis for proving weak cancellation.  Note that it is
also unknown whether real rank zero implies weak cancellation.
\endexample

\bigskip

Lawrence G. Brown           

Department of Mathematics  

Purdue University         

West Lafayette           

Indiana 47907-2067, USA 

lgb \@ math.purdue.edu

\enddocument